\documentclass[12pt]{article}
\usepackage{amssymb,amsmath,graphicx,mathtools}
\makeatletter

\@addtoreset{equation}{section}
\def\section{\@startsection {section}{1}{\z@}{-2.5ex plus -1ex minus
 -.2ex}{1.3ex plus .2ex}{\large\bf}}
\def\subsection{\@startsection{subsection}{2}{\z@}{-2.25ex plus%
 -1ex minus -.2ex}{0.5ex plus .2ex}{\bf}}

\advance \voffset by -0.8in
\advance \hoffset by -0.6in
\def\calP{{\mathcal P}}

\def\calR{{\mathcal R}}

\def\bu{{\mbox{\boldmath $u$}}}

\def\bpm{\begin{pmatrix}}
\def\epm{\end{pmatrix}}

\newcommand{\cm}{\mathfrak{m}}  
\newcommand{\cg}{\mathfrak{g}}

\newcommand{\CC}{\mathbb{C}}

\newcommand{\isom}{{\cong}}
\newcommand{\id}{{\rm id}}
\def\bee{\begin{equation}}
\def\eee{\end{equation}}
\def\bea{\begin{eqnarray}}
\def\eea{\end{eqnarray}}

\newtheorem{theorem}{Theorem}[section]

\newtheorem{lemma}[theorem]{Lemma}
\newtheorem{corollary}[theorem]{Corollary}

\def\rlbicross{{\triangleright\!\!\!\blacktriangleleft}}
\def\lrbicross{{\blacktriangleright\!\!\!\triangleleft}}

\newcommand{\C}{\mbox{${\mathbb C}$}}
\newcommand{\R}{\mbox{${\mathbb R}$}}

\newcommand{\CR}{\hbox{{$\mathcal R$}}}

\newcommand{\eps}{{\epsilon}}
\newcommand{\tens}{\mathop{\otimes}}
\newcommand{\la}{{\triangleright}}
\newcommand{\ra}{{\triangleleft}}

\newcommand{\Ad}{{\rm Ad}}
\newcommand{\ad}{{\rm ad}}
\newcommand{\<}{\langle}
\renewcommand{\>}{\rangle}

\def\rcross{{\triangleright\!\!\!<}}
\def\lcross{{>\!\!\!\triangleleft}}
\def\rcocross{{\blacktriangleright\!\!<}}

\def\rlbicross{{\triangleright\!\!\!\blacktriangleleft}}
\def\lrbicross{{\blacktriangleright\!\!\!\triangleleft}}
\def\dcross{{\bowtie}}
\def\codcross{{\blacktriangleright\!\!\blacktriangleleft}}
\def\rbiprod{{\cdot\kern-.33em\triangleright\!\!\!<}}
\def\lbiprod{{>\!\!\!\triangleleft\kern-.33em\cdot}}
\renewcommand{\o}{{}_{\scriptscriptstyle(1)}}
\renewcommand{\t}{{}_{\scriptscriptstyle(2)}}
\newcommand{\thr}{{}_{\scriptscriptstyle(3)}}
\newcommand{\fo}{{}_{\scriptscriptstyle(4)}}
\newcommand{\fiv}{{}_{\scriptscriptstyle(5)}}

\newcommand{\so}{{}^{\scriptscriptstyle[1]}}
\newcommand{\st}{{}^{\scriptscriptstyle[2]}}

\begin{document}
\parskip 4pt
\parindent 8pt

\title{Quasitriangular structure and twisting of the 2+1 bicrossproduct model}
\author{S. Majid and P. K. Osei}
\maketitle

\begin{abstract} We show that the bicrossproduct model $C[SU_2^*]\lrbicross U(su_2)$ quantum Poincare group in 2+1 dimensions acting on the quantum spacetime $[x_i,t]=\imath\lambda x_i$ is related by a Drinfeld and module-algebra twist to the quantum double $U(su_2)\rcross C[SU_2]$ acting on the quantum spacetime $[x_\mu,x_\nu]=\imath\lambda\eps_{\mu\nu\rho}x_\rho$. We obtain this twist by taking a scaling limit as $q\to 1$ of the $q$-deformed version of the above where it corresponds to a previous theory of $q$-deformed Wick rotation from $q$-Euclidean to $q$-Minkowski space. We also recover the twist result at the Lie bialgebra level. 
\end{abstract}

\section{Introduction and Motivation}

We will show that two standard and well-known quantum spacetime models, namely the `spin model' related to quantum gravity without cosmological constant and the angular momentum algebra $[x_\mu,x_\nu]=\imath\lambda\eps_{\mu\nu\rho}x_\rho$ as spacetime\cite{tHof,BatMa,FreMa,FreLev}, and the 2+1 version of the Majid-Ruegg model\cite{MaRue,AmeCam} with spacetime $[x_i,t]=\imath\lambda x_i$ are in some sense equivalent in the form of a module-algebra (or Drinfeld-type) twist. It was explained at the $q$-deformed level in \cite{MaSch} that these two models are  quantum Born reciprocal or `semidual' aspects of  3d quantum gravity and that at the $q$-deformed level (i.e. {\em with} cosmological constant) they are also related by twisting and hence in some sense self-dual up to equivalence. Quantum Born reciprocity interchanges the cosmological and Planck scales for a fixed value of $q$-deformation parameter and  the quantum double $D(U_q(su_2))=U_q(so_{1,3})$ with the quantum group $U_q(su_2)^{\text{cop}}\lrbicross U_q(su_2)=U_q(so_4)$ at the level of isometry quantum group, so these are twisting equivalent, a result first introduced in \cite{Ma:euc} as `quantum Wick rotation', see \cite{Ma:book}. Our surprising new result is that by working out the structures in great detail and carefully taking the $q\to 1$ limit while at the same time scaling the generators, i.e in a contraction limit, the result survives in the form a module algebra twist between the above two quantum spacetimes and their Poincar\'e quantum groups $D(U(su_2))$ and  $C[SU_2^*]\lrbicross U(su_2)$ respectively. The role of $D(U(su_2))$ in particular for constructing the states of 3d quantum gravity with point sources is well established and we refer to \cite{FreLev,MaSch} for an introduction. That a scaling `contraction' limit of $U_q(so_4)$ gives a quantum Euclidean group was first pointed out in \cite{Cel} and this is presumably isomorphic to $C[SU_2^*]\lrbicross U(su_2)$ in the same way as the 4d quantum Poincar\'e quantum group proposed in \cite{Luk} by contraction of $U_q(so_{2,3})$ was shown in \cite{MaRue} to be a bicrossproduct $C[\R\ltimes\R^3]\lrbicross U(so_{1,3})$. 

Our scaling limit  result is striking because the two quantum spacetimes models appear very different and have always been treated as such; one quantum spacetime is the enveloping algebra of a simple Lie algebra and the other of a solvable one. One Poincar\'e quantum group is quasitriangular while bicrossproducts are not usually quasitriangular, although the 3d one in \cite{Cel} is, a result which in our version is now explained by twisting as this preserves quasitriangularity.  Moreover, whereas the quasitriangular structure of the double exists formally, it does not take a usual algebraic form as the exponential of generators, whereas our universal R-matrix on the bicrossproduct does and this implies such a form also for the double by twisting. Similarly, when quantum spacetimes are related by a module algebra twist then their covariant noncommutative differential geometry is related by twisting\cite{MaOec,BegMa:twi} and hence that must also be the case here: For the spin model the smallest covariant calculus is known to be 4D \cite{BatMa} and for the standard bicrossproduct models it is known to be one dimension higher than classical \cite{Sit}, so again 4D but now this is explained by our twisting result.  Similarly, the construction of particle state representations of $D(U(su_2))$ by the Wigner little group method in \cite{MaSch} should have a parallel on the bicrossproduct model side via twisting. Such possible applications will be considered elsewhere. 

The paper starts in Section~2 with some general Hopf algebra constructions which underly the quantum Wick rotation\cite{Ma:euc} and semidualisation in \cite{Ma:pla,MaSch} but which were not given so explicitly before. We carefully specialise these to $U_q(su_2)$, again giving all constructions in explicit detail in Section~3.1-~3.5 . These exact formulae then  allow us in Section~3.6 to take the $q\to 1$ limit with suitably scaled generators. This is a rather tricky process due to $1/(1-q^{-2})$ singularities but we remarkably do obtain finite results, which we then verify explicitly, see Corollary~\ref{coromain}.  Section~4 rounds off the paper with the Poisson-Lie or semiclassical level version of our results in line with \cite{OS1} and mainly as a further check of our calculations (notably, we show that we recover the expected Lie bialgebra double $r$-matrices). Our results relate to a different Lie bialgebra contraction than \cite{OS2} but the latter may emerge as a different limit of our results. Another direction for future work is that $U_q(su_2)$ as quantum spacetime is a unit hyperboloid in $q$-Minkowski space and as such its  constant-time slices give the 2-parameter Podles spheres \cite{Ma:qfuzzy}, all of which may have a parallel on the bicrossproduct model side of the twisting.

\section{Explicit Hopf algebra isomorphisms}

This section brings together two different contexts in the book \cite{Ma:book}. The first, about semidualisation, was explained in \cite{MaSch} in the present context of 3d quantum gravity while the second about twisting was explained in \cite{Ma:euc} in the context of quantum Wick rotation. It was also outlined in \cite{MaSch} how to bring these together but now we need to work out the underlying isomorphisms rather explicitly, which is not easy from the literature. 

We use the conventions for Hopf algebras in \cite{Ma:book} namely a Hopf algebra or `quantum group' $H$ is both an algebra and a coalgebra, with `coproduct' $\Delta:H\to H\tens H$ which is an algebra homomorphism. There is also a counit $\eps:H\to k$ if we work over $k$ and an `antipode' $S:H\to H$ defined by $(Sh\o)h\t=h\o S h\t=\eps(h)$ for all $h\in H$ and notation $\Delta h=h\o\tens h\t$. We shall refer to a {\em covariant system} $(H,A)$ meaning a Hopf algebra $H$ acting on an algebra $A$ as a module algebra, i.e., in the left handed case,
\[ h\la (ab)= (h\o\la a)(h\t\la b),\quad h\la 1=\eps(h).\]
where $\la$ is a left action. There is then a left cross-product algebra $A\lcross H$. We refer to \cite{Ma:book} for details. We denote by $H^*$ a suitable dual Hopf algebra with dual pairing given by a non degenerate bilinear map $ \langle \,,\,\rangle$ and $H^{\text{cop}}, H^{\text{op}}$ denote taking the flipped coproduct or flipped product. As an easy exercise, if $H$ acts covariantly on $A$ from the right then  
\bea \label{flipaction} h\la a= a\ra S^{-1}h\eea is a left action of $H$ on $A^{\text{op}}$ as another covariant system. 

\subsection{Semidualisation and the quantum double}
\label{semidualisation}

(i) A {\em double crossproduct} Hopf algebra $H_1\dcross H_2$ can be thought of as a Hopf algebra $H$ which factorizes into two sub-Hopf algebras built on $H_1\tens H_2$ as a vector space. By factorisation, we mean a map $H_1\tens H_2\rightarrow H$ as an isomorphism of linear spaces.  One can then naturally  extract the actions 
$
\la: H_2\tens H_1\rightarrow H_1 \quad \mbox{and} \quad \ra: H_2\tens H_1\rightarrow H_2
$
of each Hopf algebra on the vector space of the other defined by $(1\tens a).(h\tens 1)= a\o\la h\o\tens a\t\ra h\t$
for the product  viewed on $H_1\tens H_2$ obeying some further compatibility properties (one says that one has a matched pair of interacting Hopf algebras). Conversely given such data one can reconstruct the algebra of $H_1\dcross H_2$ from these actions as a double (both left and right) cross product.  The coproduct of $H_1\dcross H_2$ is the tensor one given by the coproduct of each factor and there is  a canonical right action of this Hopf algebra on the vector space of $H_2$
which respects the coalgebra structure of $H_2$ and thus provides in a canonical way a covariant left action  of $H_1\dcross H_2$ on $H_2^* $ as an algebra. Here $H_1$ acts on $H_2^*$ by dualising the above right action $\ra$ on $H_2$, and $H_2$ acts on $H_2^*$ by the coregular action $a\la\phi=\phi\o\<a,\phi\t\>$. Hence we have a covariant system $(H_1\dcross H_2, H_2^*)$ and an associated cross product $H_2^*\lcross (H_1\dcross H_2)$. Further details are in \cite{Ma:book} and earlier works by the first author.

(ii) The  {\em semidual} of this picture  associated to the same matched pair data was introduced by the first author, see \cite{Ma:book} for details, and is  constructed by dualising the data involving $H_2$ to give a bicrossproduct Hopf algebra  $H_2^*\lrbicross H_1$ which then acts covariantly on $H_2$ from the right as an algebra as the semidual covariant system $(H_2^*\lrbicross H_1,H_2)$. The remarkable fact is that  $(H_2^*\lrbicross H_1)\rcross H_2=H_2^*\lcross (H_1\dcross H_2)$ as algebras, i.e. the combined system is the same actual algebra but its interpretation is different in that the role of spacetime coordinates $H_2$  and momentum $H_2^*$ cordinates in the first case is reversed in the other, with rotations $H_1$ the same. This is the $B$-model semidualisation referred to in \cite{MaSch}. There is equally well an $A$-model semidualisation where we dualise $H_1$ to obtain $ H_2\rlbicross H_1^*$ acting on the left on $H_1$ while $H_1\dcross H_2$ acts on the right on  $H_1^*$ and the two covariant systems again have the same cross products, $H_1\lcross(H_2\rlbicross H_1^*)= (H_1\dcross H_2)\rcross H_1^*$. These ideas go back to the first  author as a new foundation (`quantum born reciprocity') proposed for quantum gravity namely that one can swap position and momentum generators in the algebraic structure \cite{Ma:pla}. 

We will particularly need details of the $B$-model which were not provided explicitly in \cite{MaSch}. Starting with a matched pair $H_1,H_2$ acting on each other  the left action $
\la: H_1\tens H_2^*\rightarrow H_2^* 
$ of $H_1$ on $H_2^*$  and a right coaction $
 \Delta_R: H_1\to  H_1\tens  H_2^*
$ of $H_2^*$ on $H_1$ are define by
\[
(h\la\phi)(a):=\phi(a\ra h), \quad \phi\in H^*_2,\quad a\in H_2 \quad h\in H_1
\]
\[ h^0\<h^1 ,a \> =a\la h, \quad h\in H_1,\quad a\in H_2,\quad \Delta_R h=h^0\tens h^1\in H_1\tens H_2^*.\]
These define the bicrosspropduct $H_2^*\lrbicross H_1$ by a left handed cross product $H_2^*\lcross H_1$ as an algebra and a right handed cross coproduct $H_2^*\rcocross H_1$ as  coalgebra: 
\begin{align} (\phi \tens h)(\psi \tens g)=& \phi (h_{\o} \la \psi )\tens h_{\t} g,\quad \quad h\in H_1, \quad \phi,\psi\in H_2^* \\
\Delta(\phi\tens h)=&(\phi\o \tens h^0\o)\tens (\phi\t h\o^1\tens h\t)
\end{align}
The canonical right action of  $H_2^*\lrbicross H_1$ on $H_2$ is 
\bee
a\ra(\phi\tens h)=a\t\ra h \<\phi, a\o\>,\quad\forall h\in H_1,\quad a\in H_2,\quad \phi\in H_2^*.
\eee
Note that $H_2^*\tens 1$ and $1\tens H_1$ appear as subalgebras with cross relations 
 \begin{align*}
h\psi=&(1\tens h)(\psi  \tens 1)=h\o\la \psi \tens h\t
=(h\o\la \psi\tens 1)(1 \tens h\t)=(h\o\la \psi ) h\t 
\end{align*}
where we identify $h=1\tens h$  and $\psi=\psi\tens 1.$

(iii) We now apply the above construction to the specific case of the Drinfeld quantum double $D(H)=H\dcross H^{*\text{op}}$ due to \cite{Dri} and viewed as an example of a double crossproduct from work of the first author, see \cite{Ma:book} for details. Here the right action of $H$ on $H^{*\text{op}}$ and the left action of $H^{*\text{op}}$ on $H$ are given respectively by 
 \bee
 a\ra h= a\t\<h\o,a\o\>\<Sh\t,a\thr\>, \;\;  a\la h=h\t\<h\o,a\o\>\<Sh\thr,a\t\>,\; h\in H,\; a\in H^{*\text{op}}.
 \eee
The double cross product  $H\dcross H^{*\text{op}}$ then comes out as 
\bee \label{doubleproduct}
 (h\tens a).(g\tens b)=hg\t \tens b a\t \<g\o,a\o\>\<Sg\thr,a\thr\>,\quad h,g \in H, \; a,b \in H^*,\eee
 with the tensor product coproduct.  This Hopf algebra  canonically acts on $(H^{*\text{op}})^*=H^{\text{cop}}$  from the left as an algebra.  The action is
 \bee (h\tens a)\la \phi=\<\phi\o,a\> h\la  \phi\t,\quad \phi\in H^{\text{cop}} \label{DHactH}\eee
 in terms of the coproduct of $H$ and the action $\la$ in (\ref{MH(co)act}).
 
 Semidualising, the  left 
 action of $H$ on $H^{\text{cop}}$ already referred to and the right coaction of $H^{\text{cop}}$ on $H$ are respectively 
 \bee
 h\la \phi = h\o\phi S h\t=\Ad_h(\phi), \quad \Delta_R h=h\t\tens h\o Sh\thr \label{MH(co)act} 
 \eee
 so that the product  of $H^{\text{cop}}\lrbicross H$ is therefore
 \bee \label{bicrossproduct}
 (\phi\tens h)(\psi\tens g)=\phi h\o\la\psi\tens h\t g, \quad h\in H,\; \;\phi,\psi\in H^{\text{cop}}
 \eee
 as the standard cross product $H{}_\Ad\lcross H$ and the coproduct is 
  \bee\label{bicrosscoproduct}
  \Delta(\phi\tens h)=\phi\t\tens h\t\tens\phi\o h\o Sh\thr\tens h\fo
  \eee
 in terms of coproducts of $H$. This Hopf algebra acts covariantly on   $H^{*\text{op}}$ from the right according to 
  \bee\label{mirroraction}
  a\ra (\phi\tens h)=\<\phi h\o,a\o\>a\t\<S h\t, a\thr\>,\; h\in H,\; \phi\in H^{\text{cop}},\;a\in H^{*}.
  \eee
  Using \eqref{flipaction}, and correctly using the inverse antipode of the bicrossproduct determined by the coproduct  (\ref{bicrosscoproduct}) gives the covariant  left action of the bicrossproduct quantum group on $H^*$ as
  \bee\label{mirrorlaction}
   (\phi\tens h)\la a=\<Sh\o S\phi ,a\o\>a\t\< h\t, a\thr\>.
  \eee
In summary, the  semidual of the left covariant system $(D(H), H^{\text{cop}})$ is the right covariant system $(H^{\text{cop}}\lrbicross H, H^{*\text{op}})$, which is $(H^{\text{cop}}\lrbicross H, H^*)$ as a left covariant system with action (\ref{mirrorlaction}). This is essentially as in \cite{Ma:book} where we denoted $H^{\text{cop}}\lrbicross H=M(H)$ the `mirror product' but now in our current conventions and, critically, keeping track of algebras on which our Hopf algebras act.

  (iv) Finally, we observe as a right-left flipped version of \cite[Prop. 6.2.9]{Ma:book} that there is a Hopf algebra isomorphishm  
 \bee \label{theta1}
  \theta_1 :H^{\text{cop}}\tens H\rightarrow H^{\text{cop}}\lrbicross H,\; \theta_1(\phi\tens h)= \phi Sh\o\tens h\t ,\quad \theta_1^{-1}( \phi \tens h)= \phi h\o\tens h\t 
 \eee
 under which the right action of  $H^{\text{cop}}\lrbicross H$ on $H^{*\text{op}}$ by (\ref{mirroraction}) is isomorphic to a right action of $H^{\text{cop}}\tens H$  on $H^{*\text{op}}$  by 
 \begin{align*}
 a\ra (\phi\tens h)=a\ra_{H^{\text{cop}}\lrbicross H} \theta_1(\phi\tens h)
 =a\t\<\phi,a\o\>\<Sh,a\thr\>.
 \end{align*}
and by observation (\ref{flipaction}) this is equivalent to  $H^{\text{cop}}\tens H$ acting on the left on $H^*$ by 
  \begin{align}
  (\phi\tens h)\la a= a\ra S^{-1}(\phi\tens h)
  = a\ra(S\phi\tens S^{-1}h)
  =a\t\<S\phi,a\o\>\<h,a\thr\>  \label{HcopHaction}
  \end{align}
and this is also $\theta_1(\phi\tens h)\la a$ acting by (\ref{mirrorlaction}). 

In summary, the semidual of the left covariant system $(D(H), H)$ acting by (\ref{DHactH}) is isomorphic to the left covariant system $(H^{\text{cop}}\tens H, H^*)$ acting by (\ref{HcopHaction}). This action is equivalent to a left action of $H$ and a right action of another copy of $H$ it is $H^*$ with a natural $H$-bimodule structure
afforded by the coproduct (the Hopf algebra version of left and right derivatives on $H^*$).

 \subsection{Twisting of module algebras and quantum Wick rotation}
(i) We recall following Drinfeld that a quasitriangular Hopf algebra is a pair $(H,\calR),$ where $H$ is a  Hopf algebra and $\calR$ is an invertible element of $ H\tens H$ satisfying
\begin{align*}
\left(\Delta \otimes \id \right)(\calR) =& \calR_{13}\calR_{23},\quad \left(\id \otimes \Delta \right)(\calR) = \calR_{13}\calR_{12} \\
 \Delta^{\text{cop}}(h) =& \calR (\Delta h) \calR^{-1},\quad h\in H.
\end{align*}
In this case $\calR$ obeys
 \bea
 (\epsilon\tens\id)(\calR)=(\id\tens\epsilon)(\calR)=1,\nonumber \\
 (S\tens\id)\calR=\calR^{-1},\quad (\id\tens S)\calR^{-1} =\calR,\nonumber \\
 \calR_{12}\calR_{13}\calR_{23}=\calR_{23}\calR_{13}\calR_{12},\label{qybe}
 \eea
 where we write $\calR=\calR\so\tens \calR\st$ with the notation that \[\calR_{ij}=1\tens...\tens\calR\so\tens 1\tens ...\tens\calR\st\tens...\tens 1\] is the element of $H\tens H...\tens H$ which is $\calR$ in the ith and jth factors.
The identity (\ref{qybe}) is known as the {\em quantum Yang-Baxter Equation} (QYBE) and on account of this $\CR$ is also called a Universal R-matrix. If there is a $\star$-structure on $H$ the R-matrix is said to be real if $\CR^{\star\tens\star}=\CR_{21}$ and anti-real if $\CR^{\star\tens\star}=\CR^{-1}$.

Next, an element $\chi\in H\tens H$ for any Hopf algebra $H$ is called a twisting 2-cocycle \cite{Ma:book} if 
\[ \chi_{12}(\Delta\tens\id)\chi=\chi_{23}(\id\tens\Delta)\chi,\quad (\eps\tens\id)\chi=1\]
and in this case there is a new Hopf algebra $H_\chi$ with the same algebra and \cite{Ma:book},
\bea \label{Htwisted}
\Delta_\chi h = \chi(\Delta h)\chi^{-1}\quad \calR_\chi = \chi_{21}\calR {\chi}^{-1}, \quad S_\chi h = U(Sh)U^{-1}\;\;\forall h\in H_\chi,
\eea
where $U=\cdot(\id\tens S)\chi$ is invertible. Moreover, if $H$ acts covariantly on $A$ from the left then $H_\chi$ acts covariantly on a new algebra $A_\chi$ with product 
\bee \label{twistedmalgebra}
 a\cdot_\chi b=\cdot(\chi^{-1}\la (a\tens b)).
\eee
Also, if 
 $H$ is a Hopf $\star$-algebra over $\C$  and $\chi$ is real in the sense $(S \tens S)(\chi^{\star \tens \star}) = \chi_{21}$, then the twisted $\star$-structure
\bee \label{twistedstar}
\star_\chi = (S^{-1}U)(( \;)^{\star})S^{-1}U^{-1} 
\eee
 turns $H_\chi$ into a Hopf $\star$-algebra as well, see \cite{Ma:book}. 
This  cocycle twisting theory was introduced by the first author in \cite{GurMa,Ma:euc} and other works from that era (Drinfeld did not consider $2$-cocycles or module algebra twists but rather conjugation by general elements $\chi$ in the category of quasi-Hopf algebras). Clearly, if $H$ is quasitriangular and we take 2-cocycle $\chi =\calR$, then  $H_\chi=H^{\text{cop}}$.

(ii) Following \cite{Ma:euc}, we similarly see that
$H^{\text{cop}}\tens H$ acting on $H^*$ by (\ref{HcopHaction}) twists via $\chi_1=\CR_{13}^{-1}$ to $H\tens H$ acting on a new algebra, which we will denote $H^{\boxdot}$, with product 
\[ a\boxdot b = \cdot\calR\la(a\tens b)=((\CR\so\tens 1)\la a)((\CR\st\tens 1)\la b)=\calR(a\o\tens b\o)a\t b\t\] 
Thus the covariant system $(H^{\text{cop}}\tens H, H^*)$ at the end of Section~2.1 twists to $(H\tens H, H^\boxdot)$ by $\chi_1$. 

Moreover, the further twist of $H \tens H$ by the 2-cocycle $\chi_2=\calR_{23}^{-1}$ gives a Hopf algebra which we will denote  $H \codcross_\calR H$ (it is technically a double cross coproduct) acting covariantly on an algebra $\underline{H^*}$ with 
\begin{align}\nonumber
a \underline{\cdot} b=&\boxdot \left(  \calR_{23}\la (a\tens b) \right)\\\nonumber
=& a\o\<\calR\so,a\t\>\boxdot b\t \<S\calR\st,b\o\> \\\nonumber
=&a\t b\thr \calR(a\o\tens b\t)\calR(a\t\tens Sb\o)\\ \nonumber
=&a\t b\thr \calR(S a\o\tens S b\t) \calR(a\thr\tens S b\o) \\\nonumber
  =&a\t b\t \calR(Sa\o a\thr\tens S b\o)
\end{align}
for all $a,b\in H^*$, where we view $\CR$ by evaluation as a map on $H^*{}^{\tens 2}$ and use the axioms of $\CR$ in dual form. This product makes $\underline{H^*}$ with its unchanged coproduct into a braided-Hopf algebra as part of the theory of transmutation so the result in \cite{Ma:euc} was that this can be seen as a twist (namely by $\chi=\chi_2\chi_1=\CR_{23}^{-1}\CR_{13}^{-1}=(\Delta\tens\id)\CR^{-1}$). 

(iii) Moreover it is known \cite[Thm~7.3.5]{Ma:book} that there is a Hopf algebra map 
\bee \theta_2: D(H)\to H\codcross_{\CR} H,\quad \theta_2(h\tens a)=h\o \CR^{-}\st\tens h\t \CR\so\<a,\CR^{-}\so \CR\st\>
\eee
according to inclusions $i=\Delta$ and $j(a)=(\id\tens \id \tens a)(\calR^{-1}_{31}\calR_{23})$ of $H, H^{*\text{op}}$ in $H\codcross_{\CR} H$. Note also that
the latter has at least a couple of interesting quasitriangular structures built from $\CR$ namely,
 \bee
\CR_D=\calR_{41}^{-1}\calR_{31}^{-1}\calR_{24}\calR_{23}=(\chi_2)_{21}\CR_{13}^{-1}\CR_{24}\chi_2^{-1},\;\;    \calR_L=\calR_{41}^{-1}\calR_{13}\calR_{24}\calR_{23}=(\chi_2)_{21}\CR_{31}\CR_{24}\chi_2^{-1}
 \eee
 with $\CR_D$ the image under $\theta_2$ of the canonical quasitriangular structure of $D(H)$ (at least if $H$ is finite dimensional so that the latter is defined). In the factorisable case the map $\theta_2$ is an isomorphism of Hopf algebras, where `factorisable' means $Q=\CR_{21}\CR$ viewed by evaluation as a map $Q:H^*\to H$ by $Q(a)=\<a,Q\so\>Q\st$ is an isomorphism. This holds formally for the standard quantum groups associated to semisimple Lie algebras.  
 
 Pulling back under $\theta_2$ we compute using the Hopf algebra and quasitriangularity axioms that  $D(H)=H\dcross H^{*\text{op}}$ acts covariantly on $\underline{H^*}$ by
 \bea\label{DHactH*} (h\tens a)\la b=\theta_2(h\tens a)\la_{H^{\text{cop}}\tens H}b=b\thr \<h, (Sb\t)b\fo\>\CR(a\o,b\o)\CR(b\fiv,a\t)\eea
 This is also the action on of $\theta_1\theta_2(h\tens a)\in H^{\text{cop}}\lrbicross H$ on $b$ for the action (\ref{mirrorlaction}). 
 
 \begin{lemma} \label{isolemma} $Q:\underline{H^*}\to H$ is a map of covariant algebras, intertwing the action of $D(H)$ in (\ref{DHactH*}) with its action on $H$ in (\ref{DHactH}).
 \end{lemma}
 {\bf Proof} It is known \cite[Prop~7.4.3]{Ma:book} that $Q$ is a homomorphism of braided-Hopf algebras where $\underline{H^*}$ is as above with unchanged coproduct and $H$ has an unchanged product and modified coproduct $\underline\Delta$. In particular, it maps the algebras and $Q(b\t)\tens \<h,(Sb\o)b\thr\>=h\la Q(b)$ as it intertwines the left action given by evaluating with the right adjoint coaction, with the left adjoint action of $H$.  Hence
 \begin{eqnarray*} 
 Q((h\tens a)\la b)&=&\CR(a\o,b\o)\CR(b\thr,a\t) h\la Q(b\t)\\
 &=&\<Q(b\o),a\o\>\CR((Sb\t)b\fo,a\t)h\la Q(b\t)\\
 &=&\<Q(b\o)\CR\st,a\>(h\CR\so)\la Q(b\t)\\
 &=&=\<Q(b)_{\underline{(1)}}\CR\st,a\>(h\CR\so)\la Q(b)_{\underline{(2)}}\\
 &=& \<Q(b)\o,a\>h\la Q(b)\t=(h\tens a)\la Q(b)
 \end{eqnarray*}
 where we used that $Q(b\o)\tens Q(b\t)=\underline{\Delta}Q(b)=Q(b)\o S\CR\st\tens \CR\so\la Q(b)\t$ and indicated the braided coproduct by the underlining the numerical suffices. The 2nd equality is easily proven by breaking down the 3rd expressions in terms of parings of $H$ with $H^*$ and using the quasitrangular and Hopf algebra pairing axioms. $\square$

Putting all the above together, we arrive at our main result:

\begin{theorem}\label{maintheorem} If $H$ is factorisable then the covariant system $(D(H), H)$ in Section~2.1 viewed via $Q$ as a covariant system $(D(H),\underline{H^*})$ is isomorphic to a twisting of its semidual $(H^{\text{cop}}\lrbicross H,H^*)$. Here we twist by $(\theta_1\tens\theta_1)(\chi)=\CR_{23}^{-1}$ and the isomorphism is given by $\theta=\theta_1\theta_2:D(H)\to H^{\text{cop}}\lrbicross H$, where
\[ \theta(h\tens a)=h\o Q^{-}\st Sh\t\tens h\thr \CR\so\<a, Q^{-}\so \CR\st\>\]
Moreover, $H^{\text{cop}}\lrbicross H$ has two quasitriangular structures given by $\theta_1\tens\theta_1$  of $\CR_{13}^{-1}\CR_{24}$ and $\CR_{31}\CR_{24}$.
\end{theorem}
{\bf Proof} We combine the results above together with an straightforward computation for $\theta$. We recognise $\chi=\chi_2\chi_1=\CR_{23}^{-1}\CR_{13}^{-1}=(\Delta\tens\id)\CR=(\theta_1^{-1}\tens\theta_1^{-1})\CR_{23}^{-1}$ so under $\theta_1$ this maps over to $\CR_{23}^{-1}\in (H^{\text{cop}}\lrbicross H)^{\tens 2}$. Because the action of the double on the vector space $H^*$ in Lemma~\ref{isolemma} agrees via $\theta$ with the action (\ref{mirrorlaction}), it means that $Q$ at the algebra level with the transmuted product and $\theta$ at the quantum symmetry level together form an isomorphism of the covariant systems as stated. $\square$

We can also compute the quasitriangular structures in $H^{\text{cop}}\lrbicross H$ explicitly in terms of $\CR$'s using the axioms of a quasitriangular structure as 
 \begin{eqnarray*} 
\CR_{B_D}=(\theta_1\tens \theta _1) ( \CR_{13}^{-1}\CR_{24})&=& \CR_{13}^{-1}(S\tens \id \tens S\tens \id)(\Delta \tens \Delta ) \CR         \\
&=& \CR_{13}^{-1}(S\tens\id \tens S\tens \id)\CR_{14}\CR_{24}\CR_{13}\CR_{23}\\
&=&(\id\tens S^{-1}\tens\id\tens\id)( \CR_{13}^{-1}\CR_{23}\CR_{13}\CR_{14}^{-1}\CR^{-1}_{24})\\
\CR_{B_L}=(\theta_1\tens \theta _1) ( \CR_{31}\CR_{24})&=& \CR_{31}(S\tens \id \tens S\tens \id)(\Delta \tens \Delta ) \CR         \\
&=& \CR_{31}(S\tens\id \tens S\tens \id)\CR_{14}\CR_{24}\CR_{13}\CR_{23}\\
&=&(\id\tens S^{-1}\tens\id\tens\id)( \CR_{31}\CR_{23}\CR_{13}\CR_{14}^{-1}\CR^{-1}_{24})
 \end{eqnarray*}
 where all expressions are reduced to tensor products of $H$.

\section{Computations for $H=U_q(su_2)$  and $q\to 1$ scaling limit }

Here we obtain the main result, starting with explicit formulae in the $q$-deformed case. 

\subsection{The Hopf algebra $U_q(su_2)=C_q[SU_2^*]$}
We recall that the Hopf algebra $U_q(su_2)$ is defined over formal power series $\CC[[t]]$ with generators  $H$, $X_{\pm}$, where $q=e^{\frac t 2}$, say. The relations are defined by
\begin{equation}\label{uqsu2algebra}
[H, H]=0,\quad [H,X_\pm]=\pm2X_\pm,\quad [X_+,X_-]=\frac{q^H-q^{-H}}{q-q^{-1}}.
\end{equation}
The coproduct, counit and antipode are given by
\begin{align}\nonumber \label{uqsu2coalgant}
\Delta H= H\tens 1+1 \tens H \quad\quad& \Delta(X_\pm)=q^{- \frac H 2}\tens X_\pm +X_\pm\tens q^{ \frac H 2}, \\
\epsilon(H )=0 \quad\quad& \epsilon(X_\pm)=0, \\ \nonumber
S( H )=-H ,\quad\quad& S(X_\pm)=-q^{\pm 1}X_\pm,
\end{align}
For $q$ real, the $\star$-structure takes the form 
$
H^\star=H\quad X_\pm^\star= X_\mp.
$
The Hopf algebra $U_q(su_2)$ is called the $q$-deformation of the universal enveloping algebra $U(su_2)$.
It is quasitriangular with real-type universal R-matrix 
\begin{equation}\label{Uqsu2Rmatrix}
\calR=q^{\frac{ H \tens H}{2}}e^{(1-q^{-2})\,q^{ \frac H 2}X_+\tens  q^{- \frac H 2}X_-}_{q^{-2}},
\end{equation}
where $e^z_{q^{-2}}$ is the $q$-exponential
$
e^z_{q^{-2}}=\sum_{k=0}^\infty\frac{z^k}{[k;q^{-2}]!},
$
with $[k;q^{-2}]=\frac{1-q^{-2k}}{1-q^{-2}}$ and $ [k;q^{-2}]!=[k;q^{-2}][k-1;q^{-2}]...[1;q^{-2}].$  This means that $e^z_{q^{-2}}e^{-z}_{q^{2}}=1$ and that if $AB=q^2BA$, then $e^{A+B}_{q^{-2}}=e^A_{q^{-2}}e^B_{q^{-2}}.$

Next, unusually, we write $U_q(su_2)$ as $C_q[SU_2^*]$ where latter has
$\alpha,\beta,\gamma,\delta$ generators of $B_q[SU_2]$ related via  the map $Q$.
The  coordinate algebra $B_q[M_2]$  is the space of  $2\times 2$ braided Hermitian matrices \cite{Ma:exa,Ma:book}, 
or $q$-Minkowski space, with generators $\bu=\left(\begin{array}{cc}
\alpha & \beta \\
\gamma & \delta  \\
\end{array}\right)$
satisfying the relations 
\begin{align}\nonumber\label{Bqrelations}
 \beta\alpha=q^2\alpha\beta,\quad \gamma\alpha=q^{-2}\alpha\gamma,&\quad \delta\alpha=\alpha\delta,\\
\left[\beta,\gamma\right]=(1-q^{-2})\alpha(\delta-\alpha),\quad& [\delta,\beta]=(1-q^{-2})\alpha\beta,\quad [\gamma,\delta]=(1-q^{-2})\gamma\alpha,
\end{align}
and real form $ \begin{pmatrix}\alpha&\beta\\ \gamma&\delta\end{pmatrix}^\star=\begin{pmatrix}\alpha&\gamma\\ \beta&\delta\end{pmatrix}.$ Its
quotient by the braided-determinant relation
$\det(\bu)=\alpha\delta-q^2\gamma\beta=1$ gives  the braided group $B_q[SU_2]$ or $q$-hyperboloid.  When $q\ne 1$ this algebra with $\alpha^{-1}$ adjoined provides a version $U_q(su_2)$ via the map $Q$ 
 the `quantum Killing form'\cite{Ma:book} as 
\begin{equation}\label{qkillingform}
Q \begin{pmatrix}\alpha &\beta\\ \gamma &\delta\end{pmatrix}= 
\begin{pmatrix}q^H & q^{-\frac{1}{ 2}}(q-q^{-1}) q^{H\over 2}X_-\\ q^{-\frac{1}{ 2}}(q-q^{-1}) X_+q^{H\over 2}& q^{-H}+q^{-1}(q-q^{-1})^2X_+X_-\end{pmatrix}\end{equation}
which we regard as an identification. If we assume $\alpha$ is invertible then the element $\delta$ is determined by the braided-determinant relation and  not regarded as a generator. This map can also be viewed  as essentially an isomorphism between the braided enveloping algebra $BU_q(su_2)$ (which has the same algebra as $U_q(su_2)$) and its dual which is the braided function algebra $B_q[SU_2]$. 
Here, the {\em unbraided} coproduct of $B_q[SU_2]$ as inherited from that of  $U_q(su_2)$ is 
\begin{align}\nonumber\label{bgcoalgant}
\Delta \alpha &=\alpha\tens \alpha \quad \Delta\beta=1\tens \beta+\beta \tens\alpha\quad \Delta \gamma =1\tens \gamma+\gamma\tens \alpha, \\
S\alpha&=\alpha^{-1}\quad S\beta=-q^{-2}\alpha^{-1} \beta \quad S\gamma=-\gamma \alpha^{-1},\quad \epsilon(\alpha^{\pm1})=1,\quad \epsilon(\beta)=\epsilon(\gamma)=0.
\end{align}
The $R$-matrix becomes
\begin{equation}
\calR=q^{\frac{ H \tens H}{2}} e^{(1-q^{-2})^{-1}\,\gamma \tens  \alpha^{-1}\beta }_{q^{-2}}; \quad \alpha=q^H.
\end{equation}
We denote $U_q(su_2)$ in the form of the algebra of $B_q[SU_2]$ with $\alpha$ invertible and the  coproduct in \eqref{bgcoalgant} as the Hopf algebra $C_q[SU_2^*]$.

\subsection{The Hopf algebra $C_q[SU_2]=U_q(su_2^*)$}

The well-known Hopf algebra $C_q[SU_2]$ is the dual of  $U_q(su_2)$ and can be viewed as the quantum deformation of the algebra of functions of $SU(2)$. A set of generators for $C_q[SU_2]$ is given by the matrix elements $t^i{}_j:U_q(su_2)\rightarrow \C$ in the defining representation of $U_q(su_2)$ where 
\bee \label{Hpairing}
\< h, t^i{}_j\>=\rho(h)^i_j,\quad h\in U_q(su_2),\quad t^i{}_j\in C_q[SU_2]\eee
and $\rho$  is in the spin-$\frac 1 2$ representation. As usual we write
$
t^i{}_j=\left(\begin{array}{cc}
a & b \\
c & d  \\
\end{array}\right)
$
which we recall have the relations 
\begin{align}\label{uqsu2dualalgebra}\nonumber
 ba=&qab,\quad bc=cb,\quad bd=q^{-1}db,\\
 ca=&qac,\quad cd=q^{-1}dc,\quad da=ad+(q-q^{-1})bc.
 \end{align}
The coproduct, counit and antipode are given by
\begin{align}\label{Cqcoalgebraant}
 \Delta a=&a\tens a+b\tens c,\; \Delta b=a\tens b+b\tens d,\; \Delta c=c\tens a+d\tens c,\;\; \Delta d=c\tens b+d\tens d, \nonumber\\
 \epsilon a=&\epsilon d=1,\quad \epsilon b=\epsilon c=0,\quad
Sa=d,\quad Sd=a,\quad  Sb=-qb,\; \;Sc=-q^{-1}c,
 \end{align}  
and the  real form  by $a^\star=d, b^\star=-q^{-1}c$ for $q$ real.
The duality pairing takes the form
\begin{equation}\label{qpairing}
\<q^{\pm \frac H 2}, a \>=q^{\pm \frac1 2},
\quad
\<q^{\pm \frac H 2}, d \>=q^{\mp \frac1 2},
\quad
\<X_+, b \>=1,
\quad
\<X_-, c \>=1.
\end{equation} 
Applying a representation \eqref{Hpairing} to one half of the $R$-matrix  leads to the definition
of the well-known $L$-matrices
\bee
(L^+)^i_j=\calR\so\rho( \calR\st),\quad (L^-)^i_j=\rho(\calR^{-}\so) \calR^{-}\st,
\eee
where
\[
L^+= \left(\begin{array}{cc}
q^{\frac{H}{2}} &  0\\
 q^{\frac 1 2} \mu X_+ &  q^{-\frac{H}{2}}  \\
\end{array}\right)
\quad
L^-= \left(\begin{array}{cc}
q^{-\frac{H}{2}} &  - q^{\frac 3 2}\mu X_- \\
0&  q^{\frac{H}{2}}  \\
\end{array}\right),\quad \mu=1-q^{-2}.
\]
We also, unusually, write $C_q[SU_2]$ with new generators $z,x_\pm$ defined by
\begin{equation}
\left(\begin{array}{cc}
a & b \\
c & d  \\
\end{array}\right)
=
\left(\begin{array}{cc}
q^z & q^{\frac{1}{2}}\mu x_- \\
 q^{\frac{3}{2}}\mu x_+ & q^{-z} (1+q \mu^2
x_+x_- )  \\
\end{array}\right).
\end{equation}
 If we assume $a$ is invertible, the element $d$ is not regarded as a generator as it is fixed by the $q$-determinant relation $\det_q({\bf{t}})=ad-q^{-1}bc=1$.
The algebra then takes the form
\begin{equation}\label{qfuzzyrelations}
[x_\pm,z]= x_\pm,\quad  \quad [x_+,x_-]=0.
\end{equation} 
The coproduct, counit and antipode can then be translated as
\begin{align} \label{qfuzzycoalgant}
\Delta (q^z)&=q^z \tens q^z +q^2\mu^2 x_-\tens x_+ ,\quad 
\Delta(x_-)=q^z \tens x_-+x_- \tens q^{-z} +q \mu^2 x_-\tens q^{-z}x_+x_-, \nonumber\\
\Delta(x_+)&=x_+\tens q^z  + q^{-z} \tens x_+ +q \mu^2 x_+x_-q^{-z} \tens x_+, \quad
\epsilon(z)=0,\quad \epsilon(x_\pm)=0, \nonumber\\
 S(q^z)&=q^{-z}(1+q \mu^2x_+x_-),\quad S(x_\pm)=-q^{\mp 1}x_\pm.
\end{align}
We denote $C_q[SU_2]$ with $a$ invertible as the Hopf algebra $U_q(su_2^*)$. The corresponding $*$-structure on   $U_q(su_2^*)$ is given by 
 $x^\star_-=-x_+$,  $(q^z)^\star=q^{-z} (1+q \mu^2x_+x_- )$.

\subsection{The quantum double covariant system ($D(U_q(su_2)),\, U_q(su_2)$)}
The quantum double of $U_q(su_2)$ is  the double cross product Hopf  algebra $D(U_q(su_2))=U_q(su_2)\dcross C_q[SU_2]^{\text{op}}$, with
 algebra structure given by (\ref{uqsu2algebra}), the opposite algebra to (\ref{uqsu2dualalgebra}) together with cross relations obtained from \eqref{doubleproduct} as
\begin{align}\label{qdoublerelations}
[q^{\frac{H}{2}},  a]&=0,\quad q^{\frac{H}{2}} b=q^{- 1}bq^{\frac{H}{2}}, \quad q^{\frac{H}{2}} c=q c q^{\frac{H}{2}}, \quad [q^{\frac{H}{2}},  d]=0,\ \nonumber\\
X_-a&=q^{-1}a X_-+bq^{\frac H 2}, \quad [X_-,b]=0,\quad  [X_-,c]=q(q^{\frac H 2}d-q^{-\frac H 2}a),\;\; dX_-=q^{-1}X_-d+q^{\frac H 2}b,  \nonumber\\
aX_+&=q X_+a+q^{-\frac H 2}c,\quad [X_+,c]=0, \quad [X_+,b]=q^{-1}(q^{\frac H 2}a-q^{-\frac H 2}d),\;\; X_+d=qdX_++cq^{\frac H 2}.
 \end{align}
The coproduct, counit and antipode are given by \eqref{uqsu2coalgant} for the Hopf subalgebra $U_q(su_2)$ and  the coproduct, the counit and inverse of the antipode   in   \eqref{Cqcoalgebraant} for the Hopf subalgebra $C_q[SU_2]^{\text{op}}$. The quantum double $D(U_q(su_2))$ canonically acts on $U_q(su_2)$ from the left with (\ref{uqsu2algebra}) as algebra, resulting in the covariant system ($D(U_q(su_2)),\, U_q(su_2)$). The left covariance action is given by \eqref{DHactH} as
\begin{equation}\label{qdoubleaction}
\begin{array}{rl}
H \la H &=0, \quad H\la X_\pm=\pm 2 X_\pm,\quad X_\pm \la H=-2q^{\pm1}q^{-\frac H 2}X_\pm,\quad X_\pm \la X_\pm=(q^{\pm2}-q^{\pm1})q^{-\frac H 2}X^2_\pm ,\\[1.0ex]
X_\pm \la X_\mp&=q^{-\frac H 2}(X_\pm X_\mp-q^{\pm 1}X_\mp X_\pm),\quad a\, \la H =1+ H, \quad a\, \la q^{\frac{H}{ 2}}=q^{\frac1 2} q^{\frac H 2}, \quad a \,\la X_\pm=q^{-\frac 1 2}X_\pm, \\[1.0ex]
b\, \la H&=0, \quad b\, \la q^{\frac{H}{ 2}}= 0,\quad b\,\la X_+=q^{\frac H 2},\quad b\,\la X_-=0,\quad c\, \la H =0, \quad c\, \la q^{\frac{H}{ 2}}=0,  \\[1.0ex]
c\, \la X_+&=0,\quad c\, \la X_- =q^{\frac H 2}, \quad d\, \la H =-1 +H, \quad d\, \la q^{\frac{H}{ 2}}=q^{-\frac 1 2}q^{\frac{H}{ 2}}, \quad d\, \la X_\pm=q^{\frac 1 2}X_\pm. 
\end{array} \end{equation}
This is the standard q-deformed quantum double system. This $q\ne 1$ corresponds to a cosmological constant in 3d  quantum gravity.

\subsection{The bicrossproduct covariant system  $(C_q[SU^*_2]^{\text{cop}}\lrbicross U_q(su_2),\,U_q(su^*_2))$ }

Here we use the alternative description of one of the $U_q(su_2)$ as $C_q[SU_2^*]$ and of $C_q[SU_2]$ as $U_q(su_2^*)$ as explained above. From (\ref{MH(co)act}), we obtain the left action of $U_q(su_2)$ on $C_q[SU^*_2]$  as
\begin{align*} 
H\,\la \left(\begin{array}{cc}
\alpha & \beta \\
\gamma & \delta  \\
\end{array}\right)&=\left(\begin{array}{cc}
 0 & -2\beta \\
2\gamma &0  \\
\end{array}\right),\quad
q^{\pm\frac H 2}\,\la \left(\begin{array}{cc}
\alpha & \beta \\
\gamma & \delta  \\
\end{array}\right)=\left(\begin{array}{cc}
\alpha & q^{\mp1}\beta \\
q^{\pm1}\gamma & \delta  \\
\end{array}\right),\nonumber\\
X_+\,\la \left(\begin{array}{cc}
\alpha & \beta \\
\gamma & \delta  \\
\end{array}\right)&
=\left(\begin{array}{cc}
-q^{\frac 3 2}\gamma & -q^{\frac 1 2}(\delta-\alpha) \\
0 &q^{-\frac 1 2}\gamma   \\
\end{array}\right), \;\;
X_- \,\la \left(\begin{array}{cc}
\alpha & \beta \\
\gamma & \delta  \\
\end{array}\right)=\left(\begin{array}{cc}
q^{\frac{1}{ 2}}\beta & 0 \\
q^{-\frac 1 2}(\delta-\alpha) & -q^{-\frac 3 2}\beta   \\
\end{array}\right). 
\end{align*}
Writing $\alpha \equiv\alpha\tens 1$,  $\beta \equiv\beta\tens 1$, $\gamma \equiv\gamma\tens 1$, in the subalgebra  $C_q[SU^*_2]^{\text{cop}}\tens 1$ and $H\equiv1\tens H $,  $X_\pm\equiv1\tens X_\pm $,   in $1\tens U_q[su_2]$ of
the `mirror product' $M(U_q(su_2))=C_q[SU^*_2]^{\text{cop}}  \lrbicross  U_q(su_2),$ we obtain its cross relations from \eqref{bicrossproduct} as
\begin{equation}\label{qbicrossrelations}
\begin{array}{rl}
[H,\alpha]&=[H,\delta]=0,\quad [H,\beta]=-2\beta,\quad [H,\gamma]=2\gamma ,\\[1.0ex]
[X_+,\alpha]&=-q^{\frac 3 2}\gamma q^{\frac H 2},\;\; X_+\beta=q \beta X_+-q^{\frac 1 2}(\delta-\alpha)q^{\frac H 2},\;\; X_+\gamma=q^{-1}\gamma X_+,\; [X_+,\delta]=q^{-\frac 1 2}\gamma  q^{\frac H 2},\\[1.0ex]
[X_-,\alpha]&=q^{\frac1 2}\beta q^{\frac H 2},\;\; X_-\beta=q\beta X_-, \; \;X_-\gamma=q^{-1}\gamma X_-+q^{-\frac 1 2}(\delta-\alpha) q^{\frac H 2},\;\; [X_-,\delta]=-q^{-\frac 3 2}\beta q^{\frac H 2}.\\[1.0ex]
\end{array} \end{equation}
The coproduct is given by \eqref{bicrosscoproduct} as
 \begin{equation}\label{qbicrosscoproduct}
\begin{array}{rl}
\Delta \alpha &=\alpha\tens   \alpha \quad \Delta \beta =  \beta \tens 1+ \alpha \tens \beta \quad \Delta \gamma = \gamma \tens 1+ \alpha \tens  \gamma, \quad \Delta H=1\otimes H+H\otimes 1,\\[1.0ex]
\Delta X_+&=q^{-\frac H 2} \tens X_++X_+\tens \alpha^{-1}q^{\frac H 2}+q^{-\frac 1 2}\mu^{-1}(q^{\frac H 2} -q^{-\frac H 2}) \tens \gamma \alpha^{-1}q^{\frac H 2}, \\
 \Delta X_-&=q^{-\frac H 2} \tens X_-+X_-\tens \alpha^{-1}q^{\frac H 2}+ q^{-\frac 3 2}\mu^{-1}(q^{\frac H 2} -q^{-\frac H 2}) \tens  \alpha^{-1}\beta q^{\frac H 2}. 
\end{array}  \end{equation}
This Hopf algebra covariantly acts from the left on $U_q(su^*_2)$ with (\ref{qfuzzyrelations}) as  algebra giving, the covariant system $(C_q[SU^*_2]^{\text{cop}}\lrbicross U_q(su_2),\,U_q(su^*_2))$. From \eqref{mirrorlaction}, 
we obtain this  left action on $C_q[SU^*_2]$   as
\begin{align} 
H\,\la \left(\begin{array}{cc}
a& b \\
c & d  \\
\end{array}\right)&=\left(\begin{array}{cc}
 0 & -2b \\
2c &0  \\
\end{array}\right),\quad
q^{\pm\frac H 2}\,\la \left(\begin{array}{cc}
a& b \\
c & d  \\
\end{array}\right)=\left(\begin{array}{cc}
a & q^{\mp1}b \\
q^{\pm1}c & d  \\
\end{array}\right),\nonumber\\
X_+\,\la \left(\begin{array}{cc}
a& b \\
c & d  \\
\end{array}\right)&
=\left(\begin{array}{cc}
-q^{\frac 3 2}c & q^{\frac 1 2}(a- d) \\
0 &q^{-\frac 1 2}c  \\
\end{array}\right), \;\;
X_- \,\la \left(\begin{array}{cc}
a& b \\
c & d  \\
\end{array}\right)=\left(\begin{array}{cc}
q^{\frac{1}{ 2}}b & 0 \\
q^{-\frac 1 2}(d-  a) & -q^{-\frac 3 2}b  \\
\end{array}\right),\nonumber\\
\alpha \,\la \left(\begin{array}{cc}
a& b \\
c & d  \\
\end{array}\right)&=\left(\begin{array}{cc}
 q^{-1}a & q^{-1}b \\
qc &qd \\
\end{array}\right),\quad
\beta \,\la \left(\begin{array}{cc}
a& b \\
c & d  \\
\end{array}\right)=-q\mu \left(\begin{array}{cc}
0 & 0 \\
a & b  \\
\end{array}\right),\nonumber\\
\gamma \,\la \left(\begin{array}{cc}
a& b \\
c & d  \\
\end{array}\right)&
=-q^{-1}\mu \left(\begin{array}{cc}
c & d \\
0 &0\\
\end{array}\right).
\end{align}
This then translates to the left action on $U_q(su^*_2)$  given by
\begin{align} \label{bcpcovariantaction}\nonumber
H\,\la \left(\begin{array}{cc}
q^z  \\
x_+ \\
x_-  \\
\end{array}\right)&=\left(\begin{array}{cc}
0  \\
2x_+ \\
-2x_-  \\
\end{array}\right),\quad 
q^{\frac H 2}\,\la \left(\begin{array}{cc}
q^z  \\
x_+ \\
x_-  \\
\end{array}\right)=\left(\begin{array}{cc}
q^z  \\
qx_+ \\
q^{-1}x_-  \\
\end{array}\right),\\\nonumber
X_+\,\la \left(\begin{array}{cc}
q^z  \\
x_+ \\
x_-  \\
\end{array}\right)&
=\left(\begin{array}{cc}
- q^3\mu  x_+   \\
0\\
\mu^{-1}\left(q^z-q^{-z} (1+q\mu^2 x_+x_- )\right)\\
\end{array}\right),  \\\nonumber
X_- \,\la \left(\begin{array}{cc}
q^z  \\
x_+ \\
x_-  \\
\end{array}\right)&=\left(\begin{array}{cc}
q^2 \mu x_-  \\
-q^{-2}\mu^{-1}\left(q^z-q^{-z} (1+q\mu^2 x_+x_- )\right)\\
0\\
\end{array}\right),  \\ \nonumber
\alpha \,\la       \left(\begin{array}{cc}
q^z  \\
x_+ \\
x_-  \\
\end{array}\right) &=\left(\begin{array}{cc}
q^{-1}q^z  \\
qx_+ \\
q^{-1}x_-  \\
\end{array}\right),\quad \beta \,\la  \left(\begin{array}{cc}
q^z  \\
x_+ \\
x_-  \\
\end{array}\right) 
=  \left(\begin{array}{cc}
0  \\
-q^{-\frac 1 2}q^z \\
0  \\\end{array}\right),  \\ 
\gamma \,\la  \left(\begin{array}{cc}
q^z  \\
x_+ \\
x_-  \\
\end{array}\right)=& \left(\begin{array}{cc}
-q^{\frac 1 2}\mu^2x_+\\
0\\
-q^{-\frac 3 2} q^{-z} (1+q\mu^2 x_+x_- )  \\ 
\end{array}\right). 
\end{align}
This is the  `mirror product' covariant system semidual to the quantum double one. The mirror product Hopf algebra\cite{Ma:book} for generic $q\ne 1$ here is isomorphic to a tensor product so it not usually considered of interest, though it is for us.

\subsection{Twisting equivalence of the $q$-deformed covariant systems}
\label{qdeformedisom}
In this section we work out the algebra isomorphism and the twisting of the preceding two covariant systems  as established in Theorem  \ref{maintheorem}. Here the algebra isomorphism 
\[
\theta : D(U_q(su_2)) \rightarrow  C_q[SU^*_2]^{\text{cop}}
  \lrbicross  U_q(su_2)
\]
is defined in Theorem \ref{maintheorem} by
\[ \theta(h\tens {\bf{t}})=h\o Q^{-}\st Sh\t\tens h\thr \CR\so\<{\bf{t}}, Q^{-}\so \CR\st\>, \quad h\in U_q(su_2),\; {\bf{t}} \in C_q[SU_2].
\]
If  $\calR=\calR\so\tens \calR\st$, then $Q=\CR_{21}\CR=\calR\st\calR'\so\tens  \calR\so \calR'\st$, so that $Q^{-1}=\calR'^{-}\so\calR^{-}\st\tens  \calR'^{-}\st \calR^{-}\so $.
Therefore
 \begin{eqnarray*} 
Q^{-}\st\<t^i{}_k, Q^{-}\so\> &=&  \calR'^{-}\st \calR^{-}\so \<t^i{}_k, \calR'^{-}\so\calR^{-}\st  \>  \\
&=& \calR'^{-}\st \calR^{-}\so \<t^i{}_k, \calR'^{-}\so\calR^{-}\st \>  \\
&=& \calR'^{-}\st \calR^{-}\so \<t^i{}_m, \calR'^{-}\so\> \<t^m{}_k, \calR^{-}\st \>    \\
&=& \calR'\st S\calR\so \<t^i{}_m, S\calR'\so\> \<t^m{}_k, \calR\st \>    \\
&=& L^{-}{}^i{}_mSL^{+}{}^m{}_k.  
 \end{eqnarray*}
 This combination is a conjugate map $Q$ to the one we used before. Hence
 \begin{eqnarray*} 
 \theta(1\tens t^i{}_j)&=& Q^{-}\st \tens  \CR\so\<t^i{}_j, Q^{-}\so \CR\st\>=Q^{-}\st \tens  \CR\so\<t^i{}_k, Q^{-}\so\>\<t^k{}_j ,\CR\st\>\\ &=&L^{-}{}^i_mSL^{+}{}^m{}_k \tens L^{+}{}^k{}_j.
 \end{eqnarray*}
Also, for $h\tens 1\in U_q(su_2)\tens 1$, we have 
  \begin{eqnarray} 
\theta(h\tens 1)=h\o Q^{-}\st Sh\t\tens h\thr \CR\so\<1, Q^{-}\so \CR\st\> 
= h\o   Sh\t\tens h\thr =1\tens h.
 \end{eqnarray}
 In terms of our generators this means that $\theta$ identifies the $q^{H\over 2},X_\pm$ generators of the two quantum groups and 
 
 \begin{eqnarray} \label{theta}
 \theta\begin{pmatrix}a & b\\ c & d\end{pmatrix}&=&\begin{pmatrix} q^{-H}+q^3\mu^2X_-X_+ &- q^{3\over 2}\mu X_-q^{H\over 2} \\  -q^{3\over 2}\mu q^{H\over 2} X_+& q^H \end{pmatrix}\tens \begin{pmatrix}  q^{H\over 2} & 0 \\  q^{{1\over 2}}\mu X_+ & q^{-{H\over 2}} \end{pmatrix}
 \nonumber\\
 &=&\begin{pmatrix} q^2(\delta-\mu\alpha) & -q^2\beta \\ -q^2\gamma & \alpha \end{pmatrix}\tens \begin{pmatrix}  q^{H\over 2} & 0 \\  q^{{1\over 2}}\mu X_+ & q^{-{H\over 2}} \end{pmatrix}\end{eqnarray}
 where $\delta=\alpha^{-1}(1+q^2\beta\gamma) $
and in the 2nd expression we replace by  $\alpha,\beta,\gamma$ for the generators of $C_q[SU_2^*]^{\rm cop}$. One can check that this is indeed an algebra isomorphism as dictated by the theorem. 

The Drinfeld twist of the bicrossproduct $C_q[SU^*_2]^{\text{cop}} \lrbicross U_q(su_2),$ defined in Theorem \ref{maintheorem} as  $\chi_B=(\theta_1\tens\theta_1)\chi=\CR_{23}^{-1}\in C_q[SU^*_2]^{\text{cop}}\lrbicross U_q(su_2)\tens C_q[SU^*_2]^{\text{cop}}\lrbicross U_q(su_2) $
is given by \begin{align}\label{qtwist}
\chi_B&=e^{-q^{-\frac 1 2}\,1\tens q^{\frac{H}{2}}X_+ \tens  \alpha^{-1}\beta \tens 1}_{q^{2}} q^{-\frac{ 1}{2}1\tens H \tens H\tens 1}= e^{-q^{-\frac 1 2}\,K X_+ \tens  \alpha^{-1}\beta }_{q^{2}}q^{-\frac{ 1}{2} H \tens \tilde H}.
\end{align}
where $q^{\tilde H}=\alpha$ when viewed in $C_q[SU^*_2]^{\text{cop}}$ and $K=q^{H\over 2}$ in $U_q(su_2)$ and the second equality uses the identifications  $C_q[SU^*_2]^{\text{cop}}\equiv C_q[SU^*_2]^{\text{cop}}\tens 1$ and  $U_q(su_2)\equiv 1\tens U_q(su_2)$ in 
$C_q[SU^*_2]^{\text{cop}}  \lrbicross  U_q(su_2).$ 
One can check that $\chi_B(\Delta\ )\chi_B^{-1}$ where $\Delta$ is the coproduct of $C_q[SU^*_2]^{\text{cop}} \lrbicross U_q(su_2),$ gives us a coalgebra isomorphic by $\theta$ to the coalgebra of the quantum double (so $\theta$ is not just an algebra isomorphism if we take this twisted coproduct) as per the theorem. For example, 
\[ (\theta\tens\theta)\Delta d=(\theta\tens\theta)(d\tens d+c\tens b)=\alpha K^{-1}\tens \alpha K^{-1}+(-q^2\gamma K + \alpha q^{{1\over 2}}\mu X_+)\tens (-q^2\beta K^{-1})\]
\[ \chi_B(\Delta\theta(d))\chi_B^{-1}=e^{-q^{-\frac 1 2}\, K X_+ \tens  \alpha^{-1}\beta }_{q^{2}}(\alpha K^{-1}\tens \alpha K^{-1})e^{q^{-\frac 1 2}\,KX_+ \tens  \alpha^{-1}\beta }_{q^{-2}}\]
which gives the same answer on writing $A=-q^{-\frac 1 2}\, KX_+ \tens  \alpha^{-1}\beta $ and $B=\alpha K^{-1}\tens \alpha K^{-1}$ and $C=q^3 \gamma K\tens K^{-1}\beta$ so that $AB=q^2 BA+C$, $AC=CA$ using the commutation relations in $C_q[SU^*_2]^{\text{cop}} \lrbicross U_q(su_2)$.  Operators with these relations formally obey
\[ e_{q^2}^ABe_{q^{-2}}^{-A}=B+(q^2-1)BA+C,\] 
which we use. 

The map of covariant algebras provided by Lemma~\ref{isolemma}  is computed from $Q({\bf t})=\<{\bf t},Q\so\>Q\st=L^{+}SL^{-}{}$ and was already given in (\ref{qkillingform}) when one notes that $B_q[SU_2]$ and $C_q[SU_2]$ have the same coalgebra and the same generators (but different algebra relations). From this point of view $Q:U_q(su_2^*)\to U_q(su_2)$ is not an algebra map (we would have to use the transmuted or twisted product on the first algebra) and obeys
\[
Q\left(\begin{array}{cc}
q^z & q^{\frac{1}{2}}\mu x_- \\
 q^{\frac{3}{2}}\mu x_+ & q^{-z} (1+q \mu^2
x_+x_- )  \\
\end{array}\right) = 
\begin{pmatrix}q^H & q^{\frac{1}{ 2}}\mu q^{H\over 2}X_-\\ q^{\frac{1}{ 2}}\mu X_+q^{H\over 2}& q^{-H}+\mu^2X_+X_-\end{pmatrix}.
\]
Here $a^n=a\underline\cdot a \underline\cdot a\cdots$ ($n$ times) when one looks carefully at the transmuted product $\underline\cdot$ on the generator $a=q^z$, which implies 
\bee \label{qcovmap} Q(z)=H,\quad Q(x_-)=q^{H\over 2}X_-,\quad Q(x_+)=q^{-1}X_+q^{H\over 2}.  \eee
at the level of $U_q(su_2^*)$ generators.

We also obtain two quasitriangular structures for the bicrossproduct $C_q[SU^*_2]^{\text{cop}} \lrbicross  U_q(su_2)$ defined in Theorem \ref{maintheorem}. Then we find expressions for  $$\calR_{B_D},\CR_{B_L}\in C_q[SU^*_2]^{\text{cop}}  \lrbicross U_q(su_2)\tens C_q[SU^*_2]^{\text{cop}}  \lrbicross U_q(su_2)$$  as follows:
From Theorem \ref{maintheorem}, we have $\CR_{B_D}=(\id\tens S^{-1}\tens\id\tens\id)( \CR_{13}^{-1}\CR_{23}\CR_{13}\CR_{14}^{-1}\CR^{-1}_{24})$ and similarly for $\CR_L$ with $\CR_{31}$ in place of $\CR_{13}^{-1}$. Writing $\CR$ for $U_q(su_2)$, this is
\begin{align} \label{RBD}
\CR_{B_D}&=(\id\tens S^{-1}\tens\id\tens\id)\big(q^{{1\over 2}1\tens H\tens H\tens 1}e_{q^2}^{-\mu KX_+\tens q^H\tens K^{-1}X_-\tens 1}e_{q^{-2}}^{\mu q^H \tens KX_+\tens K^{-1}X_-\tens 1} \nonumber\\
&e_{q^{-2}}^{\mu K X_+\tens 1\tens K^{-1}X_-\tens 1}
e_{q^2}^{-\mu KX_+\tens 1\tens1\tens K^{-1}X_-}e_{q^2}^{-\mu q^H\tens K X_+\tens 1\tens K^{-1}X_-}q^{-{1\over 2}(H\tens1\tens 1\tens H+1\tens H\tens 1\tens H)}\big)\nonumber\\
&=(\id\tens S^{-1}\tens\id\tens\id)\big(q^{{1\over 2}1\tens H\tens \tilde{H}  \tens 1}e_{q^{-2}}^{q^{-\frac 1 2} \alpha \tens KX_+\tens \alpha^{-1}\beta \tens 1}e_{q^2}^{-\frac 1 \mu \gamma \tens q^H\tens \alpha^{-1}\beta\tens 1} \nonumber\\
&e_{q^{-2}}^{\frac 1 \mu \gamma \tens 1\tens \alpha^{-1}\beta \tens 1}
e_{q^2}^{-q^{\frac 1 2}\gamma \tens 1\tens1\tens K^{-1}X_-}e_{q^2}^{-\mu \alpha \tens K X_+\tens 1\tens K^{-1}X_-}q^{-{1\over 2}(\tilde{H}\tens1\tens 1\tens H+1\tens H\tens 1\tens H)}\big)
\end{align}
where  $K=q^{H\over 2}$ and $\alpha=q^{\tilde{H}}$ viewed in $C_q[SU^*_2]^{\text{cop}}$ .  Here, we have written $\calR_{B_D}$ as an element of $U_q(su_2)  \lrbicross U_q(su_2)\tens U_q(su_2)  \lrbicross U_q(su_2)$ in the first equality.
In the second equality, we view the first and the third legs in $C_q[SU^*_2]^{\text{cop}}$ using the map $Q$ in \eqref{qkillingform} and used the fact that $e_{q^{-2}}^{q^{-\frac 1 2} \alpha \tens KX_+\tens \alpha^{-1}\beta \tens 1}$ commutes with $e_{q^2}^{-\frac 1 \mu \gamma \tens q^H\tens \alpha^{-1}\beta\tens 1}$.
Note that  $S^{-1}$  reverses order, resulting in more complicated expressions if we apply this. Similarly,  Theorem \ref{maintheorem} gives
\begin{align} \label{RBL}
\CR_{B_L}&=(\id\tens S^{-1}\tens\id\tens\id)\big(q^{{1\over 2}H\tens 1\tens H\tens 1}e_{q^{-2}}^{\mu K^{-1}X_-\tens 1 \tens KX_+\tens 1}q^{{1\over 2}1\tens H\tens H\tens 1}e_{q^{-2}}^{\mu 1 \tens KX_+\tens K^{-1}X_-\tens 1} \nonumber\\
&q^{{1\over 2}H\tens 1\tens H\tens 1} e_{q^{-2}}^{\mu K X_+\tens 1\tens K^{-1}X_-\tens 1}
e_{q^2}^{-\mu KX_+\tens 1\tens1\tens K^{-1}X_-}q^{-{1\over 2}H\tens 1\tens 1\tens H}\nonumber\\
&e_{q^2}^{-\mu 1\tens K X_+\tens 1\tens K^{-1}X_-}q^{-{1\over 2}(1\tens H\tens 1\tens H)}\big)\nonumber\\
&=(\id\tens S^{-1}\tens\id\tens\id)\big(q^{{1\over 2}\tilde H\tens 1\tens \tilde H\tens 1}e_{q^{-2}}^{\frac 1 \mu \alpha^{-1}\beta \tens 1 \tens \gamma\tens 1}q^{{1\over 2}1\tens H\tens \tilde H\tens 1}e_{q^{-2}}^{q^{-\frac 1 2} 1 \tens KX_+\tens \alpha^{-1}\beta \tens 1} \nonumber\\
&q^{{1\over 2}\tilde H\tens 1\tens \tilde H\tens 1} e_{q^{-2}}^{\frac 1\mu \gamma\tens 1\tens \alpha^{-1}\beta \tens 1}
e_{q^2}^{-q^{\frac 1 2} \gamma\tens 1\tens1\tens K^{-1}X_-}q^{-{1\over 2}\tilde H\tens 1\tens 1\tens H}\nonumber\\
&e_{q^2}^{-\mu 1\tens K X_+\tens 1\tens K^{-1}X_-}q^{-{1\over 2}(1\tens H\tens 1\tens H)}\big).
\end{align}

 \subsection{Limiting twist between the spin model and the bicrossproduct model}
 \label{limitingtwist}
 We are now in position to consider the degenerations of the two covariant systems by scaling the various generators appropriately to recover the results of Theorem \ref{maintheorem}  in the limit $q\rightarrow 1$. 
 
 (i) For the quantum double $D(U_q(su_2))=U_q(su_2)\dcross  C_q[SU_2]^{\text{op}}$,  the $U_q(su_2)$ part has  no problem with the limit $q\rightarrow 1$  and so $U_q(su_2)\mapsto U(su_2)$,  with the standard Lie brackets
 \begin{equation}\label{su2algebra}
[H, H]=0,\quad [H,X_\pm]=\pm2X_\pm,\quad [X_+,X_-]=H,
\end{equation}
and cocommutative coalgebra
\begin{align}\nonumber \label{su2coalgant}
\Delta H&= H\tens 1+1 \tens H \quad\quad \Delta(X_\pm)=1\tens X_\pm +X_\pm\tens 1, \\
\epsilon(H )&=0 \quad \epsilon(X_\pm)=0,\quad\quad
S( H )=-H ,\quad S(X_\pm)=-X_\pm,
\end{align}
Similarly,  $C_q[SU_2]^{\text{op}} \mapsto C(SU_2)$, the commutative algebra of functions on  $SU_2$ as $q\rightarrow 1$. 
The cross relations in the limit can easily be extracted from \eqref{qdoublerelations} as
\begin{align}\label{doublerelations}\nonumber
[H,  a]=0,\quad [H, b]=- b, &\quad [H,c]=2c, \quad [H,  d]=0,\ \\\nonumber
[X_-,a]=b, \quad [X_-,b]=0,&\quad  [X_-,c]=d-a,\;\; [X_-,d]=-b,  \\\nonumber
[X_+,a]=-c,\quad [X_+,c]=0, &\quad [X_+,b]=a-d),\;\; [X_+,d]=c,\\
 \end{align}
 and the coproduct is the tensor product one.  The quantum  group $U_q(su_2)$ is also the quantum spacetime algebra for the covariant system ($D(U_q(su_2)),\, U_q(su_2)$). In the limit
 the covariant action of $D(U(su_2))=U(su_2)\rcross C(SU_2)$ on $U(su_2)$ is given by the $q\to 1$ limit of \eqref{doubleaction} as
 \begin{equation}\label{doubleaction}
\begin{array}{rl}
H \la H &=0, \quad H\la X_\pm=\pm 2 X_\pm,\quad X_\pm \la H=-2X_\pm,\quad X_\pm \la X_\pm=0 ,\\[1.0ex]
X_\pm \la X_\mp&=H,\quad a\, \la H =1+ H,  \quad a \,\la X_\pm=X_\pm,\quad b\, \la H=0, \quad b\,\la X_+= 1,\quad \\[1.0ex]
 b\,\la X_-&=0,\quad c\, \la H =0, \quad c\, \la X_+=0,\quad 
c\, \la X_- =1, \quad d\, \la H =-1 +H, \quad  d\, \la X_\pm=X_\pm. 
\end{array} \end{equation}
This limit action was first computed in \cite{BatMa} and to match standard conventions, we choose generators 
\bee  H=2  J_0,\quad X_\pm=   J_1\pm \imath J_2,\quad t^i_j=\left(\begin{array}{cc}
e^{\lambda \calP_3}+\imath \frac {\lambda}{ 2}\calP_0 & \imath \frac {\lambda}{ 2}(\calP_1-\imath\calP_2) \\
\imath \frac {\lambda}{ 2}(\calP_1+\imath\calP_2) & e^{\lambda \calP_3}-\imath \frac {\lambda}{ 2}\calP_0  \\
\end{array}\right),
\eee
where $\lambda\in \R$ is a deformation parameter and $\calP_3$  not regarded as a generator but is determined by the $\det(t)=1$ condition. This gives the algebra as
\bee \label{Dconventiongen}
[J_a,J_b]=\imath \eps_{abc}J_c,\quad [\calP_a,J_b]=\imath \eps_{abc}\calP_c,\quad [\calP_a,\calP_b]=0,
\eee
and the coproducts turns out to be
\bee
\Delta J_a=J_a\tens 1+1\tens J_a,\quad \Delta \calP_a=\calP_a\tens 1+1\tens \calP_a-\frac \lambda 2  \eps_{abc}\calP_b\tens \calP_c.
\eee
 For  the quantum spacetime algebra for the covariant system,  we  write
\bee \lambda H=2  x_0,\quad \lambda X_\pm=   x_1\pm \imath x_2
\eee
and then the limit of relations \eqref{su2algebra}  gives the spin model spacetime  algebra
\bee \label{ncspace}
[x_\mu,x_\nu]=\imath\lambda\eps_{\mu\nu\rho}x_\rho.
\eee
This is the $q\to 1$ limit $(D(U(su_2)),\, U(su_2))$ as a deformation of $U(iso(3))$ on $\R^3$\cite{BatMa}.

(ii) In the covariant system $(C_q[SU^*_2]^{\text{cop}}\lrbicross U_q(su_2),\,U_q(su^*_2))$,  we have  $C_q[SU^*_2]^{\text{cop}}\mapsto C[SU^*_2]^{\text{cop}}$ with commutative algebra and $U_q(su_2)\mapsto U_(su_2)$ with the standard algebra \eqref{su2algebra}.
 From \eqref{qbicrossrelations}, the cross relations becomes
 \begin{equation}\label{limitbicrossrelations}
\begin{array}{rl}
[H,\alpha]&=[H,\delta]=0,\quad [H,\beta]=-2\beta,\quad [H,\gamma]=2\gamma ,\\[1.0ex]
[X_+,\alpha]&=- \gamma,\;\; [X_+,\beta]= -(\delta-\alpha),\;\; [X_+,\gamma]=0,\; [X_+,\delta]=\gamma,\;\;  \\[1.0ex]
[X_-,\alpha]&=\beta,\;\; [ X_-,\beta]=0, \; \;[X_-,\gamma]=(\delta-\alpha),\;\; [X_-,\delta]=-\beta.\\[1.0ex]
\end{array} \end{equation}
The coproduct is obtained from \eqref{qbicrosscoproduct} as
 \begin{equation}\label{limitbicrosscoproduct}
\begin{array}{rl}
\Delta \alpha &=\alpha\tens   \alpha \quad \Delta \beta =  \beta \tens 1+ \alpha \tens \beta \quad \Delta \gamma = \gamma \tens 1+ \alpha \tens  \gamma, \quad \Delta H=1\otimes H+H\otimes 1\\[1.0ex]
\Delta X_+&=1 \tens X_++X_+\tens \alpha^{-1}+\frac H 2\tens \gamma \alpha^{-1}, \;\;
 \Delta X_-=1 \tens X_-+X_-\tens \alpha^{-1}+ \frac H 2 \tens  \alpha^{-1}\beta . 
\end{array}  \end{equation}
In this limit the covariant action of $C[SU^*_2]^{\text{cop}}\lrbicross U(su_2)$ on   $U(su^*_2)$  from \eqref{bcpcovariantaction} is
\begin{align}\label{BCPlimitactions}
H\la \,z&=0,\;\; H\la \,x_\pm=\pm 2x_\pm,\;\; X_\pm\la\, x_\pm=0,\;\; X_\pm\la\, z=\mp 2x_\pm,\;\;X_\pm\la\,x_\mp=\pm z \nonumber\\
 \alpha^{\pm1}\la\, z&=z\mp1,\quad \alpha\la\,x_\pm=x_\pm,\quad  \beta\la \,z=0,\;\; \beta \la \,x_+=-1,\;\; \beta \la\, x_-=0,\;\; \gamma \la \,z=0\nonumber\\
 \gamma \la \,x_+&=0,\;\; \gamma \la\, x_-=-1.
\end{align}
To match standard conventions for the bicrossproduct quantum group $C[SU^*_2]^{\text{cop}}\lrbicross U(su_2)$, we now identify new generators
\begin{eqnarray} \label{limitGEN} \nonumber
\alpha&=&e^{\lambda p_0},\quad \beta = \lambda P_+,\quad \gamma =\lambda P_-,\quad P_\pm= p_2\pm \imath  p_1,\\ 
H&=& 2M,  \quad X_+=N_2- \imath N_1,\quad X_-=N_2+ \imath N_1.
\end{eqnarray}
Then the relations \eqref{limitbicrossrelations}  become 
\begin{equation}\label{qto1bicrossrelations}
\begin{array}{rl}
\left[p_a,p_b\right] &=0, \quad [M,N_1]=N_2, \quad [M,N_2]=-N_1, \quad [N_1,N_2]=-M,\\[1.0ex]
\left[M,p_0\right]&=0 \quad \left[M, p_i\right]=\imath \epsilon_{ij}p_j,\quad \left[N_i,p_0,\right]=-\imath\epsilon_{ij}p_je^{-\lambda p_0} \\[1.0ex]
\left[N_i,p_j\right]&=\frac \imath 2\epsilon_{ij}\,e^{-\lambda p_0}\left(\frac{ e^{2\lambda p_0}-1}{\lambda}- \lambda \vec{p}\,^2 \right),\quad i,j=1,2,  \\[1.5ex]
\end{array} \end{equation}
where $\vec{p}\,^2=p_1^2+p_2^2$ and  coproduct  \eqref{limitbicrosscoproduct} gives
 \begin{equation}
\begin{array}{rl}
\Delta p_0 &=p_0\otimes1+1\otimes p_0,\quad \Delta M=1\otimes M+M\otimes 1,\\[1.0ex]
\Delta p_i &= p_i\tens 1+e^{\lambda p_0} \otimes p_i ,\\[1.0ex]
\Delta N_i&=1\otimes N_i+N_i \otimes e^{-\lambda p_0}+ \lambda M\tens p_ie^{-\lambda p_0},\quad i=1,2. \end{array}  \end{equation}
For the model spacetime $U_q(su^*_2)$, we set
 \bee
 x_0=\imath\lambda z,\quad x_1=-\imath\lambda( x_++x_-),\quad x_2=\lambda( x_+-x_-)
 \eee
  and then take the limit $q\rightarrow 1$ to get from \eqref{qfuzzyrelations}, the bicrossproduct model spacetime
\bee
[x_i,x_0]= \imath \lambda x_i.
\eee
 In terms of these standard generators \eqref{limitGEN}, the covariant actions \eqref{BCPlimitactions} can be translated as 
\begin{align}
M\la \,x_0&=0,\;\; M\la \,x_i=-\epsilon_{ij}x_j,\;\;  N_i\la\, x_0=-\imath x_i,\;\;  N_i\la\, x_j=-\imath \delta_{ij}x_0,\;\;  \nonumber\\
p_0 \la x_0&=-\imath,\quad  p_0\la\, x_i=0,\quad  p_i\la \,x_0=0,\quad p_i\la \,x_j=i\epsilon_{ij},\;\; i,j=1,2.
\end{align}
This is the bicrossproduct model covariant system $(C[SU^*_2]^{\text{cop}}\lrbicross U(su_2),U(su_2^*))$ as a quantum Poincare group in three dimensions acting on the Majid-Ruegg quantum spacetime as a 3d version of \cite{MaRue}. 

(iii) We now look at the $q\to1$ limit of the twist between the two covariant systems. We remind the reader that for handling of the cocycle and R-matrices we reduced expressions to the tensor product of the underling Hopf algebras. In effect in what follows we equip the vector space of $C[SU^*_2]^{\text{cop}}\lrbicross U(su_2)$ with two products, one is the cross product algebra as part of the bicrossproduct $\lrbicross $ construction and the other is the tensor product $\tens$ algebra. In the following, we use the convention that all exponentials are multiplied in the tensor product $\tens$ algebra, which does not impact (\ref{qto1twist}) but is important for the correct reading of (\ref{qto1RBD}). 
\begin{corollary} \label{coromain}
From the above analysis, we arrive at the degeneration limit $q\rightarrow 1$ of our result that the covariant system $(D(U(su_2))=U(su_2)\dcross  C[SU_2]^{\text{op}},\, U(su_2)$) is isomorphic to a twisting of  the covariant system $(C[SU^*_2]^{\text{cop}}\lrbicross U(su_2),\,U(su^*_2))$. 
The twist is derived from \eqref{qtwist} in the limit $q\rightarrow 1$ as 
\bee \label{qto1twist}
\chi_{B_0}= e^{- X_+ \tens \alpha^{-1}\beta}\, e^{-\frac 1 2H \tens (\alpha -1)}\in (C[SU^*_2]^{\text{cop}}\lrbicross U(su_2))^{\tens 2}.
\eee
The degeneration limit $q\rightarrow 1$  of \eqref{RBD} also provides an R-matrix for the bicrossproduct $C[SU^*_2]^{\text{cop}}\lrbicross U(su_2)$ given by  
\bee \label{qto1RBD}
\CR_{B_0}=  e^{\frac 1 2  \gamma H\tens\alpha^{-1}\beta }  : e^{- \alpha X_+\tens \alpha^{-1}\beta }:   e^{-{1\over 2} H\tens (\alpha-1) }e^{-\gamma  \tens X_-}  e^{-{1\over 2} (\alpha-1) \tens H} ,
\eee
where $$: e^{- \alpha X_+\tens \alpha^{-1}\beta }:=\sum_{n=0}^\infty \frac{ (-1)^n} {n!} \alpha^n  X_+^n\tens (\alpha^{-1}\beta)^n.$$
\end{corollary} 
{\bf Proof} We write $\alpha \equiv\alpha\tens 1$,  $\beta \equiv\beta\tens 1$, $\gamma \equiv\gamma\tens 1$, in the subalgebra  $C_q[SU^*_2]^{\text{cop}}\tens 1$ and $H\equiv1\tens H $,  $X_\pm\equiv1\tens X_\pm $,   in $1\tens U_q(su_2)$ of
$C_q[SU^*_2]^{\text{cop}}  \lrbicross  U_q(su_2),$
with $q=e^t$. Then in the limit $t\to0$, the algebra isomorphism  becomes
$ \theta : D(U(su_2)) \rightarrow  C[SU^*_2]^{\text{cop}}
  \lrbicross  U(su_2)
$
given by
 \bee \label{isomlambda}
 \theta\begin{pmatrix}a & b\\ c & d\end{pmatrix}
 =\begin{pmatrix} \delta & -\beta \\ -\gamma & \alpha \end{pmatrix},\quad \theta(h)=h, \; h\in U(su_2)
 \eee
It is easy to check that $\theta$ is indeed an algebra isomorphism. For example $[\theta(X_+),\theta(b)]=\theta(a)-\theta(d),$ etc.
Now we write  $\alpha=e^{{t\over 2}\tilde{H}}=1+{t\over 2}\tilde{H}+O(t^2)$, so that 
$$q^{-\frac{ 1}{2} H \tens \tilde H}=\sum_{n=0}^\infty \frac {(-t)^n} {n!} \frac {{H^n}\tens {\tilde{H}^n}}{4^n}=\sum_{n=0}^\infty \frac {1} {n!}  \left(-\frac{H}{2}\right)^{n}\tens (\alpha-1)^n+O(t^2)=e^{-\frac 1 2 H\tens (\alpha-1)}+O(t^2),$$ and therefore from \eqref{qtwist}, we see that $\chi_{B}\to \chi_{B_0}=e^{- X_+ \tens \alpha^{-1}\beta}\,e^{-\frac 1 2H \tens (\alpha -1)} $ as $t\to 0$. To obtain the limit for the $R$-matrix, we first note that for $X$ commuting with $H$, 
$$\lim_{t\to 0} \left(e_{q^{2}}^{-\frac X \mu e^{{t\over 2}H}}e_{q^{-2}}^{\frac X \mu }\right)=\sum_{n=0}^\infty \frac {1} {n!} \left(-\frac X 2\right)^n\prod^{n-1}_{i=1}(H-2i).$$ Then form \eqref{RBD}, the limit of  $\calR_{B_D}$
 as $q\to 1$ becomes
\begin{align} \label{limitstep1}
\calR_{B_0}=&(\id\tens S \tens\id\tens\id) \Big(e^{{1\over 2}1\tens H\tens (\alpha-1)\tens 1 } e^{ \alpha\tens X_+\tens \alpha^{-1}\beta\tens 1 }\nonumber\\
&\sum_{n=0}^\infty \frac {(-1)^n} {n!2^n}  \left(\gamma\tens1 \tens \alpha^{-1}\beta\tens 1 \right)^{n}(1\tens \prod^{n-1}_{i=1}(H-2i)\tens 1\tens 1)
\nonumber\\
&e^{-\gamma\tens 1\tens 1 \tens X_-}
e^{-{1\over 2}(\alpha-1)\tens 1\tens 1 \tens H}\Big),
\end{align}
where we used that $X=\gamma\tens1\tens\alpha^{-1}\beta\tens 1$ commutes with $1\tens H\tens 1\tens 1$ and also that $S^2=\id$ on $U(su_2)$.
Next, we also note that for any elements
\begin{align*}
(\id\tens S \tens\id\tens\id)&\big((a\tens b\tens c\tens d)\cdot(A\tens B\tens C\tens D)\big)=(\id\tens S \tens\id\tens\id)(aA\tens bB\tens  cC\tens dD)\\
&=aA\tens S(bB)\tens cC\tens dD=Aa\tens (SB)(Sb)\tens Cc\tens Dd\\
&=(\id\tens S \tens\id\tens\id)(A\tens B\tens C\tens D)\cdot (\id\tens S \tens\id\tens\id)(a\tens b\tens c\tens d)
\end{align*}
provided $dD=Dd$, since the first and third legs are in $C[SU_2^{*}]^{\text{cop}}$ which is already commutative. Here,  $\cdot$ indicates that the product is in the tensor product one of the Hopf algebra.
Using this observations, \eqref{limitstep1} becomes
\begin{align} \label{limitstep2}
\calR_{B_0}=
&(\id\tens S \tens\id\tens\id)\big(\sum_{n=0}^\infty \frac {1} {n!2^n}  \gamma^n\tens\prod^{n-1}_{i=1}(H+2i) \tens( \alpha^{-1}\beta)^n\tens 1) \big)
\nonumber\\
&(\id\tens S \tens\id\tens\id) \big(e^{{1\over 2}1\tens H\tens (\alpha-1)\tens 1 } e^{ \alpha\tens X_+\tens \alpha^{-1}\beta\tens 1 }\big)\nonumber\\
&e^{-\gamma\tens 1\tens 1 \tens X_-}e^{-{1\over 2}\left((\alpha-1)\tens 1\tens 1 \tens H \right)}
\end{align}
because the first two lines of \eqref{limitstep1} are of the form $A\tens B\tens C\tens 1$ and the last line is unchanged under the action of $(\id\tens S \tens\id\tens\id)$.
We evaluate the first two lines of the above equation  as follows: In the first line, we have 
\begin{align*} 
&(\id\tens S \tens\id\tens\id)\sum_{n=0}^\infty \frac {(-1)^n} {n!2^n}  \gamma^n\tens\prod^{n-1}_{i=1}(H-2i) \tens( \alpha^{-1}\beta)^n\tens 1
\nonumber\\
&=\sum_{n=0}^\infty \frac {1} {n!2^n}  \gamma^n\tens\prod^{n-1}_{i=1}(H+2i) \tens( \alpha^{-1}\beta)^n\tens 1
 =e^{\frac 1 2  \gamma H\tens\alpha^{-1}\beta }.
\end{align*}
The second line of \eqref{limitstep2} is  evaluated as
\begin{align*} 
&(\id\tens S \tens\id\tens\id) \big(e^{{1\over 2}1\tens H\tens (\alpha-1)\tens 1 } e^{ \alpha\tens X_+\tens \alpha^{-1}\beta\tens 1 }\big)\nonumber\\
&=e^{ -\alpha\tens X_+\tens \alpha^{-1}\beta\tens 1 }e^{-{1\over 2}1\tens H\tens (\alpha-1)\tens 1 }=\sum_{n=0}^\infty \frac{ (-1)^n} {n!} \alpha^n  X_+^n\tens (\alpha^{-1}\beta)^ne^{-{1\over 2} H\tens (\alpha-1) }\\
&=
: e^{- \alpha X_+\tens \alpha^{-1}\beta }:  e^{-{1\over 2} H\tens (\alpha-1) }.
\end{align*}

Finally, we check that  $U(su_2)$ is a module algebra twist by $\chi_{B_0}$ of $U(su_2)$, i.e. $U(su^*_2)_{\chi_{B_0}}=U(su_2)$, where the twisted product $a\cdot_{\chi_{B_0}} b$ is given by \eqref{twistedmalgebra}.
With
$ \chi^{-1}_{B_0}= e^{\frac H 2 \tens (\alpha -1)}\,e^{ X_+ \tens \alpha^{-1}\beta},
$
we have that for example
\begin{align*} 
z \cdot_ {\chi_{B_0}} x_+&=\cdot({\chi^{-1}_{B_0}}\la (z\tens x_+))\\
&=\sum\sum\frac 1 {m!}\frac 1 {n!}\left(\left(\frac H 2\right)^m\la X^n_+ \la z   \right)\left((\alpha-1)^m\la (\alpha^{-1}\beta)^n\la x_+\right)\\
&=\sum \frac 1 {m!}\left(\left(\frac H 2\right)^m\la z  \right)\left((\alpha-1)^m\la  x_+\right)+\sum \frac 1 {m!}\left(\left(\frac H 2\right)^m\la (-2x_+)   \right)\left((\alpha-1)^m\la  (-1)\right)\\
&=zx_++2x_+,
\end{align*}
and 
\begin{align*} 
x_+ \cdot_ {\chi_{B_0}} z&=\sum\sum\frac 1 {m!}\frac 1 {n!}\left(\left(\frac H 2\right)^m\la X^n_+ \la x_+   \right)\left((\alpha-1)^m\la (\alpha^{-1}\beta)^n\la z\right)\\
&=\sum \frac 1 {m!}\left(\left(\frac H 2\right)^m\la x_+  \right)\left((\alpha-1)^m\la  z \right)\\
&=x_+z-x_+.
\end{align*}
Computing all  possible combination of products, we obtain the following twisted algebra
\begin{align*}
[x_+,x_-] \cdot_ {\chi_{B_0}} &=[x_+,x_-]+z,\quad   [x_\pm,x_\pm] \cdot_ {\chi_{B_0}} =[x_\pm,x_\pm],\quad [z,z] \cdot_ {\chi_{B_0}} =[z,z] \nonumber\\
  [x_+,z] \cdot_ {\chi_{B_0}} &=[x_+,z]-3x_+,\quad  [x_-,z] \cdot_ {\chi_{B_0}} =[x_-,z]+x_-=2x_-
 \end{align*}
which on evaluating  the product in $U(su^*_2)$ gives
\[
[x_+,x_-] \cdot_ {\chi_{B_0}} =z,\quad   [x_\pm,x_\pm] \cdot_ {\chi_{B_0}} =0,\quad [z,z] \cdot_ {\chi_{B_0}} =0,\quad [x_\pm,z] \cdot_ {\chi_{B_0}} =\mp 2x_\pm. 
\]
We see that
\bee
\Phi(H)=z,\quad \Phi(X_\pm)=x_\pm
\eee
defines an isomorphism of  $U(su_2)$ with $U(su^*_2)$ after twisting.  This is manifestly the inverse of $Q:U(su_2^*)\to U(su_2)$ of covariant algebras in Lemma~\ref{isolemma} which in our case by \eqref{qcovmap} is just $Q(z)=H$, $Q(x_\pm)=X_\pm$ in the $q\to 1$ limit. By construction, the identification must be covariant but it is a useful check to see this directly. For example, 
\[ \theta(X_+)\la \Phi(H)=X_+\la z=-2x_+=\Phi(X_+\la H)\]
\[ \theta(a)\la \Phi(H)=\delta\la z=z+1=\Phi(1+H)=\Phi(a\la H)\]
\[ \theta(d)\la \Phi(X_-)=\alpha \la x_-=\Phi(X_-)=\Phi(d\la X_-).\]

\section{Semiclassical limit of results}

In this section describe our twisting result at the infinitesimal Lie bialgebra level, i.e. the classical double   as a Lie bialgebra twist of the bicross sum,  both in the case when $q$  is switched on and in the scaling limit $q\to 1$ with a parameter $\lambda$.  We begin with a  brief review of Lie bialgebras and classical $r$-matrices and   refer the reader to  \cite{Ma:book,CP} and references therein for details.

\subsection{Double and Semidual Lie Bialgebras and classical $r$-matrices}\label{sec4.1}
A Lie bialgebra in the sense of Drinfeld provides a semiclassical or infinitesimal notion of a Hopf algebra.
A Lie bialgebra $(\cg,[\;,\;],\delta)$ is a Lie algebra $(\cg,[\;,\;])$ over a field $k$ equipped with a cocommutator $\delta:\cg\mapsto \cg\otimes  \cg$  is a skew-symmetric linear map, i.e. $\delta:\cg\mapsto \wedge^2\cg$ satisfying the the coJacobi identity 
 \[
(\delta\otimes \text{id})\circ \delta(\xi)+\mbox{cyclic}=0, \quad\forall \xi\in \cg
\]
and that for all $\xi,\eta\in \cg$, 
\[  \delta([\xi,\eta])=(\ad_\xi\otimes 1+1\otimes \ad_\xi)\delta(\eta)-(\ad_\eta\otimes 1+1\otimes \ad_\eta)\delta(\xi).
 \]
There exist a Lie bialgebra  version of the quasitriangular Hopf algebra which arise naturally in the following way: Since $\delta$ is a 1-cocycle,
an element $r= r \so\tens r\st \in \cg\tens \cg$ provides a coboundary structure for the Lie bialgebra $(\cg,[\;,\;],\delta)$ if $\delta = \partial r$, i.e.
 $ \delta(\xi)=\ad_\xi(r)=[\xi\otimes1+1\otimes \xi,r]$. 
This requires that $\ad_\xi(r+r_{21})=0$ for all $\xi\in\cg$ to have $\delta$ antisymmetric. For any Lie algebra $\cg,$ we define the  map  \begin{equation}
\cg^{\otimes2}\rightarrow \cg^{\otimes3},\;\;r\mapsto[[r,r]]=[r_{12,}r_{13}]+[r_{12,}r_{23}]+[r_{13,}r_{23}].
\end{equation}
This map restricts to the map $\wedge^{2}\cg\rightarrow\wedge^{3}\cg.$
The equation \begin{equation}
\label{CYBE} [[r,r]]=0
\end{equation} is called the classical Yang-Baxter equation(CYBE) and any solution of the CYBE in $\cg\otimes \cg$ is called a classical $r$-matrix.
The classical $r$-matrix provides a quasitriangular structure for the Lie bialgebra. It is triangular if it satisfies the CYBE and $r_{21}=-r$ and called factorisable if it satisfies the CYBE and $r+r_{21}: \cg^*\rightarrow \cg$ is a linear surjection. The classical r-matrix therefore provides a natural infinitesimal version of the universal R-matrix while the factorisables case correspond to $\calR$ factorisable.

If $(\cg,[\;,\;],r)$ is a quasitriangular Lie bialgebra  and $\chi^c\in\cg\tens\cg$ obeys
\bee 
[[r,\chi^c]]+[[\chi^c,r]]+[[\chi^c,\chi^c]]=0,\quad \ad_\xi(\chi^c+\chi^c_{21})=0,\quad \forall \xi\in \cg,
\eee
then $(\cg,[\;,\;],r+\chi^c)$ is also a quasitriangular Lie bialgebra. The element $\chi$ is called a Lie bialgebra  twist and $\delta_{\chi^c}=\delta + \partial\chi^c$ is also a Lie bialgebra.

A quantised enveloping algebra roughly speaking means a Hopf algebra over $\CC[[t]]$ generated by a vector space $\cg$, with relations and coproduct of the form
\bee
\xi\eta-\eta\xi=[\xi,\eta]+O(t), \quad (\Delta-\Delta^{\text{cop}})\xi= t\delta \xi+O(t^2),
\eee
where $t$ is a formal deformation parameter.
 Further, if the Hopf algebra has a quasitriangular structure of the form
\bee \calR = 1 + t r + O(t^2),\label{Rexpantion}
\eee
then $\delta =\partial r$, and our Lie bialgebra is quasitriangular.  This interpretation is also compatible with twisting. More explicitly, if $\chi=1+tf+O(t^2)$, then from \eqref{Htwisted}, we have 
\bee \label{ctwist}
\chi^c=f_{21}-f
\eee

Classical double  and bicross sum Lie bialgebras provide a  semiclassical version of the quantum doubles and bicrossproduct quantum groups respectively, described in Section \ref{semidualisation}.
For any  finite dimensional Lie bialgebra $\cg$ with dual $\cg^*$, there is a quasitriangular Lie bialgebra, $D(\cg)$, the classical double of $\cg$ built on $\cg\oplus\cg^*$ as a vector space, with
\begin{eqnarray*}
[\xi\oplus\phi,\eta\oplus\psi] &=&([\xi,\eta] +\sum\xi_{\so}\langle\psi,\xi_{\st}\rangle-\eta_{\so}\langle\phi,\eta_{\st}\rangle) \\
&\quad &  \oplus ([\psi,\phi] +\sum\psi_{\so}\langle \psi_{\st},\xi\rangle -\phi_{\so}\langle\phi_{\st},\eta\rangle),\\
\delta(\xi\oplus\phi) &=& \sum(0\oplus\phi_{\so})\tens(0\oplus\phi_{\st})+\sum(\xi_{\so}\oplus 0)\tens(\xi_{\st}\oplus 0),\\
\quad  r&=&\sum_{a}(0\oplus f^{a})\tens(e_{a}\oplus 0).
\end{eqnarray*}
Here $\cg^{*\text{op}},\cg$,  appear as sub-Lie bialgebras, where $()^{\text{op}}$ denotes the opposite (negated) Lie bracket. The set $\lbrace e_{a}\rbrace$ is a basis of $\cg$ and $\lbrace f^{a}\rbrace$, a dual basis. Moreover, if $(\cg,r)$ is factorisable then $D(\cg)\isom (\cg\oplus\cg)_{\chi_2^c}$ where we mean twisting by a certain cocycle $\chi_2^c$, which as we saw amounts to adding $\chi_2^c$ to the $(-r_{21})\oplus r$ if we want $r_{D}$ and to $r\oplus r$ for another $r$-matrix $r_{L}$. We can further view this as a twisting of $\cg^{\text{cop}}\oplus\cg$ by a certain other cocycle $\chi_1^c$ built from $r$. 

Semidualisation can also be defined for Lie bialgebras which are double cross sums.
Given a matched pair of Lie algebras $(\cg,\cm),$ one can define the double cross sum $\cg\dcross \cm$ as the vector space $\cg \oplus \cm$. The semidual gives the the bicross sum Lie bialgebra $\cm^*\lrbicross \cg$ built on $\cm^* \oplus \cg$
with 
\begin{equation*}
[f\oplus \xi,h\oplus\eta] =(\xi\la h -\eta\la f)\oplus[\xi,\eta],
\end{equation*}
\begin{eqnarray*}
\delta(f\oplus \xi) &=& \sum_{a}(0\oplus e_{a}\la\xi)\tens(f^{a}\oplus 0)-(f^{a}\oplus 0)\tens(0\oplus e_{a}\la\xi)\\
&\quad & +\sum (f_{\so}\oplus 0)\tens(f_{\st}\oplus 0),
\end{eqnarray*}
for all $f\oplus \xi, h\oplus \eta\in \cm^{*}\lrbicross \cg$, where $\delta(f)=\sum f_{\so}\tens f_{\st}$ is the Lie colagebra given by the dualization of the Lie bracket of $\cm$. For a detailed account of these constructions, we refer to \cite{Ma:book}. In particular, the splitting from the double semidualises to a bicross-sum $\cg^{\text{cop}}\lrbicross \cg\isom \cg^{\text{cop}}\oplus\cg$ as Lie bialgebras. 

Putting these facts together gives an isomorphism of quasitriangular Lie bialgebras $\theta^c:D(\cg)\to (\cg^{\text{cop}}\lrbicross \cg)_{\chi_B^c}$, with
\bee \label{thetac} 
\theta^c(\xi)=\xi,\quad \theta^c(\phi)=-2\tilde r_+(\phi)+(\id\tens\phi)r,\quad \chi_B^c=r_{23}-r_{41},\eee
for all $\xi\in \cg$ and $\phi\in \cg^*$, where $r_+=(r+r_{21})/2$ is viewed as a map $\cg^*\to \cg$ and tilde indicates that the result is viewed in the $\cg^{\text{cop}}$ copy. As a check, one has
\[ (\theta^c\tens\theta^c)(f^a\tens e_a)=r_{24}-r_{14}-r_{41}=r_{B_D}+\chi_B^c,\quad r_{B_D}=r_{24}-r_{23}-r_{14}\]
so the canonical $ r$-matrix for the double when mapped over under the isomorphism is the  twist of $r_{B_D}$. There is also 
\[ r_{B_L}=r_{31}-r_{23}+r_{13}-r_{14}+r_{24}\]
which twists to the image under $\theta^c\tens\theta^c$ of the other $r$-matrix $r_{L}$ on $D(\cg)$. This is the Lie bialgebra version of the general theory in Section~2. We now verify everything on our examples as a check.

\subsection{Infinitesimal  limit of results in the limit $q\to 1$}
In the infinitesimal Lie bialgebra limit, the quantum double $D(U(su_2))$ becomes the  Lie bialgebra double $D(su_2)=su_2\rcross su_2^{*{\text{op}}}=su_2\rcross \R^3$. Here, the  $su_2$  parts has  its standard Lie bracket and $su^{*\text{op}}_2=\R^3$ has a commutative algebra. The relations are given by \eqref{Dconventiongen}
\bee \label{}
[J_a,J_b]=\imath \eps_{abc}J_c,\quad [\calP_a,J_b]=\imath \eps_{abc}\calP_c,\quad [\calP_a,\calP_b]=0.
\eee

The bicrossproduct quantum group $C[SU^*_2]^{\text{cop}}\lrbicross U(su_2)$ described in Section \ref{limitingtwist} in terms of the basis \eqref{limitGEN} can be viewed  as a deformation of $ U(su^{\text{cop}}_2\lrbicross  su_2 )$, where $\lambda$ is the deformation parameter.
In the semiclassical limit, $C[SU^*_2]^{\text{cop}}\lrbicross U(su_2)$ becomes the bicross sum $su^{\text{cop}}_2\lrbicross  su_2$, where $su_2$ has  its standard Lie bracket and $su^{\text{cop}}_2=\R^3$ has a commutative Lie bracket. From \eqref{qto1bicrossrelations} the relations are given by 
\begin{equation}
\begin{array}{rl}
\left[p_a,p_b\right] &=0, \quad [M,N_1]=N_2, \quad [M,N_2]=-N_1, \quad [N_1,N_2]=-M,\\[1.0ex]
\left[M,p_0\right]&=0 \quad \left[M, p_i\right]=\imath \epsilon_{ij}p_j,\quad \left[N_i,p_0,\right]=-\imath\epsilon_{ij}p_j,\quad \left[N_i,p_j\right]= \imath \epsilon_{ij}p_0 \quad i,j=1,2, \\[1.0ex]
\end{array} \end{equation}

In terms of the basis \eqref{limitGEN},
 the $R$-matrix becomes
\bee \label{RPMN}
\CR_{B_0}=  e^{  \lambda P_- M\tens \lambda P_+ e^{-\lambda p_0} }  : e^{- e^{\lambda p_0} X_+\tens \lambda P_+e^{-\lambda p_0}}:   e^{-M\tens (e^{\lambda p_0}-1) }e^{-\lambda P_- \tens X_-}  e^{- (e^{\lambda p_0}-1) \tens M}  ,
\eee
where we have kept the $X_\pm$ for simplicity.
Then the semiclassical limit of the bicrossproduct $R$-matrix \eqref{RPMN} gives a classical $r$-matrix for the bicross sum  as
\bea \label{bicrossrmatrix1}
r_{B_0}&=&   - M\tens  p_0  -  p_0  \tens M    - P_- \tens X_-    - X_+\tens  P_+ \nonumber\\
&=&  - M\tens  p_0  -  p_0  \tens M    -(p_2\tens N_2+p_1\tens N_1+N_2\tens p_2+N_1\tens p_1)\nonumber\\
&&-\imath(p_2\tens N_1-p_1\tens N_2+N_2\tens p_1-N_1\tens p_2)
\eea
It is interesting to see that if we set $M=J_0$, $J_i=N_i$ and $P_a=-p_a$, to match standard notation, we get an $r$-matrix for the bicross sum $su^{\text{cop}}_2\lrbicross  su_2$ as
\bea \label{bicrossrmatrix2}
r_{B_0}&=&   P_a\tens J_a +J_a\tens P_a -\imath (P_1\wedge J_2-P_2\wedge J_1)\nonumber\\
&=& P_a\tens J_a +J_a\tens P_a +\imath m_a\epsilon^{abc} P_b\wedge J_c,\quad {\bf {m}}^2=1,
\eea
where $\bf m$ is a unit time-like vector. 

We observe that the symmetric part of the $r$-matrix \eqref{bicrossrmatrix2} is equal to the Casimir associated to the invariant, non-degenerate symmetric bilinear form used in the Chern-Simons action \cite{AT,Witten}  and therefore suitable for constructing the Poisson structure on the classical phase space via the Fock-Rosly construction \cite{FockRosly}. This shows that the bicrosproduct  with $r$-matrix depending  on a time-like deformation vector and with complex antisymmetric part is compatible with 3d gravity via the Fock-Rosly compatibility condition. This  is different from the family of classical bicross sum $r$-matrices  associated to 3d gravity with vanishing cosmological constant obtained in \cite{MSkappa,OS1}. In the later, the $r$-matrices are real and depend on space-like deformations vectors.  See also \cite{OS2}  where a complete classification of all $r$-matrices compatible to 3d gravity with vanishing cosmological is constructed via semidualisation of  Lie bialgebras which are double cross sums.

Now, rewriting  the twist in Corollary \ref{coromain} in terms of the basis \eqref{limitGEN}, we get
\bee \label{twistPMN}
\chi_{B_0}= e^{- X_+ \tens \lambda P_+e^{-\lambda p_0}}\, e^{-M \tens (e^{\lambda p_0} -1)},
\eee
and the semiclassical limit  for the twist gives the Lie bialgebra twist by \eqref{ctwist} as 
\begin{align} \label{CtwistB0}
\chi^c_{B_0}&= X_+ \tens  P_+ +M \tens  p_0-P_+\tens X_+-p_0\tens M  \nonumber \\
&=M \tens  p_0-p_0\tens M+N_2\tens p_1+N_1\tens p_1-p_2\tens N_2-p_1\tens N_1\nonumber\\
&+\imath(N_2\tens p_1-N_1\tens p_2-p_2\tens N_1-p_1\tens N_2)
\end{align}

The above considerations leads to the semiclassical limit  of the results in Corollary~\ref{coromain} that the double Lie bialgebra $D(su_2)$ is a Lie bialgebra twisting of  bicross sum Lie bialgebra  $su_2^{\text{cop}}\lrbicross su_2$. The isomorphism \eqref{isomlambda} becomes 
$
\theta^c:D(su_2)\to su_2^{\text{cop}}\lrbicross su_2,
$ 
where
\bee \label{classicalisom}
\theta^c(J_0)=M,\;\; \theta^c(J_i)=N_i,\quad  \theta^c(\calP_a)=-2 p_a,\;\; i=1,2, \;\; a=0,1,2.
\eee
Twisting the bicrossproduct $r$-matrix by $\chi^c_{B_0}$  and using the isomorphism \eqref{classicalisom} gives the $r$-matrix for a classical double $D(su_2)$ as
\bea 
r_{D_0}=r_{B_0} +\chi^c_{B_0} &=& -2p_0\tens M-P_-\tens X_--P_+\tens X_+=  \calP_a\tens J_a
\eea
on using the identification $\theta^c$ for the last step, in agreement with the general theory in Section~\ref{sec4.1}.

\subsection{Infinitesimal  limit of the $q$-deformed  results}

If we take the semiclassical limit of Section \ref{qdeformedisom} without sending $q\to 1$, we have on the one side the Lie bialgebra double $D(su_2)=su_2\dcross su_2^{*op}$ where $su_2^{*op}=an_2$ is the Lie algebra of the  Lie group $AN_2$ of  $2\times 2$ matrices of the form
\begin{equation*}
\label{an2para}
\begin{pmatrix}
e^{\phi} & \xi+i\eta \\
0 & e^{-\phi} \\
 \end{pmatrix},\quad \phi ,\xi,\eta \in \R,
 \end{equation*}
and the notation refers the abelian and the nilpotent parts of this group. The vector space splitting expressed in the Lie double cross sum is the lie version of the Iwasawa factorisation of $SL(2,\CC)=SU_2.AN_2$. This result from \cite{Ma:mat} is the reason that the quantum double $D(U_q(su_2))$ can be regarded as the $q$-Lorentz quantum group. Our result is that this $D(su_2)$ is a twist of the double cross sum $su_2\lrbicross su_2$ as quasi-triangular Lie bialgebras. The latter is known to be isomorphic to $su_2^{\text{cop}}\oplus su_2$ recovering the (complexified) Lorentz Lie bialgebra as twist of a direct sum. This fact is essentially known\cite{Ma:book} but it is a nice check of our formulae to check this from the semiclassical limit of Section \ref{qdeformedisom}. 

To this end, the quantum double is  given by (\ref{uqsu2algebra}), (\ref{qfuzzyrelations}), with cross relations in  \eqref{qdoublerelations} written conveniently as 
\begin{align}\label{uqsu2dualalgebra2}
[q^{\frac{H}{2}},  q^z]&=0,\quad q^{\frac{H}{2}} x_\pm=q^{\pm 1}x_\pm q^{\frac{H}{2}},\quad  [X_\pm,x_\pm]=0, \nonumber\\
[X_-,x_+]&= q^{-\frac{1}{2}} \mu^{-1}  \left( q^{\frac{H}{2}} q^{-z}(1+q \mu^2 x_+x_-)- q^{-\frac{H}{2}}q^{z} \right),\;\; X_- q^z= q^{-1} q^zX_- +q^{\frac1 2}\mu x_-q^{\frac{H}{2}}, \nonumber\\
  [X_+,x_-]&=  q^{-\frac{3}{2}}\mu^{-1}  \left( q^{\frac{H}{2}} q^{-z}(1+q \mu^2 x_+x_-)- q^{-\frac{H}{2}} q^z \right),\quad q^zX_+= q X_+ q^z+q^{\frac 3 2}\mu q^{-\frac{H}{2}}x_+ .
 \end{align}
The coproduct, counit and antipode are given by \eqref{uqsu2coalgant} for  $U_q(su_2)$ and the opposite of the coproduct, the counit and inverse of the antipode   in  (\ref{qfuzzycoalgant})  for $U_q(su^*_2)^{\text{op}}$.   The 
bicrossproduct  $U_q(su_2)^{\text{cop}}\lrbicross U_q(su_2)$  is from \eqref{qbicrossrelations} and \eqref{qbicrosscoproduct}.
The twisting \eqref{qtwist} takes the form
 \bee\label{qtwist2}
\chi'_B=e^{-(1-q^{-2})\,1\tens q^{\frac H 2}X_+ \tens q^{-\frac {\tilde{H}}{ 2}} \tilde{X}_-\tens 1}_{q^{2}}  q^{-\frac 1 2  1 \tens H \tens \tilde{H}\tens 1}=e^{-(1-q^{-2})\,q^{\frac H 2}X_+ \tens q^{-\frac {\tilde{H}}{ 2}} \tilde{X}_-}_{q^{2}}  q^{-\frac{ 1}{ 2}   H \tens \tilde{H}},
\eee
 where $\tilde{H}, \tilde{X}_\pm$ are the generators for $U_q(su_2)^{\text{cop}}$.

For the semiclassical or infinitesimal regime, in \eqref{RBD} and \eqref{RBL}, we write $q=e^{\frac t 2}$ and use  \eqref{Rexpantion} to get    classical $r$-matrices for the bicross sum  as
\bee \label{rBD}
r_{B_D}= \frac 1 4  ( H \tens H-\tilde{H}\tens H- H\tens \tilde{H} ) +X_+\tens X_- -\tilde{X}_+ \tens X_- -X_+  \tens \tilde{X}_- ,
 \eee
 and
\begin{eqnarray}\label{rBL}
r_{B_L}&=& \frac 1 4  (2\tilde{H}\tens \tilde{H}-\tilde{H}\tens H-H\tens \tilde{H}+ H\tens H )\nonumber\\ 
&+& \tilde{X}_-\tens\tilde{X}_++ \tilde{X}_+\tens\tilde{X}_-+X_+\tens  X_-- X_+ \tens \tilde{X}_-  -\tilde{X}_+  \tens X_- .
 \end{eqnarray}
 It is easy to check that $r_{B_D}$ and $r_{B_L}$ satisfy the CYBE \eqref{CYBE}.
 In this semiclassical regime, the twisting \eqref{qtwist2} becomes 
\bee
\chi_B^c=\frac 1 4  (H\tens\tilde{H}- \tilde{ H} \tens H)  + X_+ \tens  \tilde{X}_--\tilde{X}_-\tens X_+.
\eee
Thus the two classical $r$-matrices for the corresponding classical double are
\bee  \label{rDD}
r_{D}=\chi_B^c+r_{B_D} = \frac 1 4  (H \tens H-2\tilde{H}\tens H )  +X_+\tens X_--\tilde{ X}_+ \tens  X_--\tilde{X}_-\tens X_+,
\eee
\begin{align}  \label{rDL}
r_{L}=\chi_B^c+r_{B_L} &= \frac 1 4  (H \tens H -2\tilde{H}\tens H  -2\tilde{H}\tens\tilde{H})\nonumber\\
&+\tilde{X}_-\tens\tilde{X}_++ \tilde{X}_+\tens\tilde{X}_-+X_+\tens  X_-- \tilde{X}_+ \tens X_- - \tilde{X}_- \tens X_+.
\end{align}

Now for the Lie algebra $su_2$ with basis $\{H,X_\pm\}$,  the standard Drinfeld-Sklyanin r-matrix is $r=\dfrac 1 4 H\tens H +X_+\tens X_-$ so that  $r_+= \dfrac 1 4 H\tens H +\dfrac 1 2(X_+\tens X_-+X_-\tens X_+)$.  We let $su_2^*=\text{span}\{\phi,\psi_\pm\}$  be the dual Lie algebra with relations
\[ [\psi_\pm, \phi] = \frac 1 2 \psi_\pm,\quad [\psi_+, \psi_-] = 0,
\]
and dual pairing  
\[\<\phi, H\>=1,\quad \<\psi_+, X_+\>=1,\quad \<\psi_-, X_-\>=1.
\]
Then from \eqref{thetac} we get
\[
\theta^c(H)=H,\; \theta^c(X_\pm)=X_\pm,\; \theta^c(\phi)=-\frac{\tilde{H}} 2 + \frac{H} 4,\; \theta^c(\psi_+)=-\tilde{X}_-,\; \theta^c(\psi_-)=-\tilde{X}_++X_+.
\]
which one can check is in agreement with semiclassicalising \eqref{theta}. Then
\begin{align*}(\theta^c\tens \theta^c)(\phi \tens H&+\psi_- \tens X_- +\psi_+ \tens X_+)\\ 
&=( \frac H 4 -\frac{\tilde{H}}{  2})\tens H +(X_+-\tilde{X}_+) \tens  X_-  -\tilde{X}_-\tens X_+=r_D\end{align*}
in agreement with the general theory in Section~\ref{sec4.1}. 

\subsection*{Acknowledgement}
PKO thanks Perimeter Institute and Fields Institute for the Fields-Perimeter Africa Postdoctoral Fellowship.
This research was partly supported by the Perimeter Institute for Theoretical Physics. Research at Perimeter Institute is supported by the Government of Canada through the Department of Innovation, Science and Economic Development Canada and by the Province of Ontario through the Ministry of Research, Innovation and Science.

\end{document}